\documentclass[a4paper,oneside, 11pt,reqno]{amsart}
\usepackage[
    a4paper,
    left=2.5cm,
    right=2.5cm,
    top=3.0cm,
    bottom=3.0cm
]{geometry}
\allowdisplaybreaks
\usepackage{amssymb,amsmath,mathtools,amsfonts,mathrsfs,calc,physics}
\usepackage{subfigure}
\usepackage{xcolor}
\usepackage[dvipsnames]{xcolor}
\usepackage{epsfig}
\usepackage[all]{xy}
\usepackage[utf8]{inputenc}
\usepackage{graphicx}
\usepackage{xcolor}
\usepackage{comment}
\usepackage{stmaryrd}
\usepackage{enumitem}
\usepackage{dynkin-diagrams}
\usepackage{hyperref} 
\usepackage{ragged2e}
\usepackage{vwcol}
\usepackage{paracol} 
\usepackage{lipsum}
\usepackage{xparse}
\usepackage[pro]{fontawesome5}
\usepackage{tcolorbox}
\usepackage[]{mdframed}
\newtcolorbox{mybox}[1]{colback=blue!3!white,colframe=blue!60!black,fonttitle=\bfseries, title=#1}
\usepackage{faktor}
\usepackage{footmisc} 
\usepackage{matlab-prettifier}
\hypersetup{colorlinks=true,linkcolor=blue,citecolor=RubineRed,urlcolor=blue}
\usepackage{anyfontsize}    
\usepackage[linesnumbered,ruled,noline]{algorithm2e} 
\SetNlSty{textnormal}{}{:}   
\SetAlgoNlRelativeSize{-1}
\SetKwComment{Comment}{\textcolor{black}{$\triangleright$ }}{\textcolor{ForestGreen}{}}
\usepackage{matlab-prettifier}

\newcommand{\mc}[1]{\mathcal{#1}}
\newcommand{\mb}[1]{\mathbb{#1}}

\newcommand{\mfrac}[2]{\scalebox{0.9}{$\dfrac{#1}{#2}$}}

\newtheorem*{claim}{Claim} 
\newenvironment{proofofclaim}%
{\begin{proof}}%
  {\end{proof}}

\newtheorem{theorem}{Theorem}[section] 
\newtheorem{lemma}[theorem]{Lemma}  
\newtheorem{corollary}[theorem]{Corollary}
\newtheorem{proposition}[theorem]{Proposition}

\theoremstyle{definition}

\newtheorem{definition}[theorem]{Definition}

\newtheoremstyle{lemmastyle}
  {3pt}
  {3pt}
  {\itshape}
  {}
  {\normalfont\scshape}%
  {.}
  {.5em}
  {}
\theoremstyle{lemmastyle}
{\begin{proof}}%
  {\end{proof}}

\title[On Nonlocal Equations under Symmetries]{On Nonlocal Equations with Symmetries}

\begin{document}

\maketitle

\begin{center}
\author{\textsc{Santiago Castellanos}\footnote{
Universidad de Santiago de Chile, Santiago de Chile.
Email: \email{santiago.castellanosg10@gmail.com}.\\
This article is based on the research conducted during the author's Ph.D. thesis \cite{Castellanos 1}.}}
\end{center}

\begin{abstract}
In this work, we establishes results of existence and regularity for solutions to nonlocal equations of the form $f(\Delta)u = U(x,u(x)) ~\text{and}~f(\Delta_{g})u=U(x,u(x))$
in the Euclidean context and Riemannian context respectively. We extend the results of \cite{Gorka 6} for the Euclidean space $\mathbb{R}^{n}$ under cylindrically symmetry and compact Riemannian manifolds $(M,g)$ under action of a closed Lie subgroup $G \leq \text{Isom}_{g}(M)$. Our approach combines spectral theory, complex interpolation and fixed-point methods.
\vspace{0.2cm}

\justify
\end{abstract}

\section{Introduction}

In this article, we study generalized nonlocal equations, both linear and 
nonlinear, of the form
\begin{equation}\label{ecuacion no local}
    f(\Delta)u = U(x, u(x)),
\end{equation}
where $f$ is an arbitrary measurable real-valued function defined on 
$\mathbb{R}$, and $U$ is a 
prescribed nonlinearity. Our central objective is to prove that, under 
appropriate hypotheses on $f$ and $U$, equations of this form admit smooth 
solutions under the influence of symmetries. The investigation of nonlocal equations of type 
\eqref{ecuacion no local} is strongly motivated by their relevance in
contemporary mathematical physics, particularly in \emph{string theory} and 
\emph{nonlocal cosmological models}\cite{Kostelecky}. Among the most paradigmatic examples is the 
generalized \emph{bosonic string equation}
\begin{equation}
    \Delta e^{-c\Delta}\phi - U(x,\phi) = 0, \qquad c > 0.
\end{equation}
Another 
central example is the \emph{open $p$-adic string equation} for the scalar tachyon 
field,
\begin{equation}
    p^{a\partial_t^2}\phi = \phi^p, \qquad a > 0,
\end{equation}
studied rigorously in \cite{Vladimirov 2}.
At the level of generality of \eqref{ecuacion no local} neither pseudo-differential operators in the 
classical sense of H\"{o}rmander \cite{Hormander} nor the theory of classical symbols or conditions are necessary.
We seek a direct meaning for the operator 
$f(\Delta)$, or its equivalent on Riemannian manifolds; to this end, we represent the nonlocal part of the equation \eqref{ecuacion no local} via the Fourier transform, functional calculus, and spectral theory, following the approach initiated by G\'{o}rka, Prado, 
and Reyes in \cite{Gorka 6}. 
Once interpreted, we proceed to solve the nonlocal equation \eqref{ecuacion no local} in different settings in the presence of symmetries, either cylindrical symmetry in Euclidean spaces or isometric group actions on compact Riemannian manifolds, on which the operator $f(\Delta) + I$ (resp., $f(\Delta_g) + I$) is well defined and acts naturally.

For the Euclidean context, the conditions imposed on $f: \mathbb{R}\to\mathbb{R}$ are the following:\\
$(P)$~~\text{The function} $s \mapsto f(-s^{2})$~\text{is nonnegative}.\\
$(E_{\beta})$ There exist real numbers $\beta, R, M > 0$ such that
$$M(1+|\xi|^{2})^{\frac{\beta}{2}} \leq f(-|\xi|^{2})~~\forall \xi~\text{such that}~ |\xi|> R.$$
Functions satisfying these conditions are called $\mathcal{G}^{\beta}$-symbols, and we denote the class of such functions by $\mathcal{G}^{\beta}$. 
In the Euclidean setting, the nonlocal operator $f(\Delta)$ is interpreted through the Fourier transform as
$
f(\Delta)u=\mathcal{F}^{-1}\bigl(f(-|\xi|^2)\mathcal{F}(u)(\xi)\bigr),
$
so that $f$ is evaluated at $-|\xi|^2\in(-\infty,0]$. 
When working with compact Riemannian manifolds $(M,g)$, the nonlocal equation we work with varies slightly with respect to \eqref{ecuacion no local}, since the nonlocal term is now $f(\Delta_{g})$ where $\Delta_{g}$ is the Laplace-Beltrami operator defined as $\Delta_{g}=-\text{div}_{g}(\nabla_{g})$, and the conditions on $f$ are formulated in terms of the discrete spectrum $(\lambda_i)_{i\ge 1}$ of $\Delta_g$.
\vspace{0.2cm}

In this article, we study the problem under the influence of symmetry in two main contexts:
\renewcommand{\theenumi}{\Roman{enumi}}
\begin{enumerate}
    \setcounter{enumi}{0}
    \item For the Euclidean case under the influence of cylindrical symmetry in the sense of \cite{Lions}, we first recall the ambient space without symmetry: for a fixed symbol $f \in \mc{G}^{\beta}$, we work with the Hilbert space $\mc{H}^{\beta}(f)$ given by all measurable real-valued functions $g$ on $\mathbb{R}^n$ whose Fourier transform $\mathcal{F}(g)$ exists and satisfies
    \begin{equation}
   \int_{\mathbb{R}^n} \left[ 1 + f\left(-|\xi|^2\right) \right]^2 \left| \mathcal{F}(g)(\xi) \right|^2 \, d\xi < \infty.
    \end{equation}
This space satisfies $\mathcal{H}^{\beta}(f) \subseteq L^2(\mathbb{R}^n)$, and allows us to give an interpretation of $f(\Delta) + I$. 
To incorporate cylindrical symmetry, we work with the spaces  
    $$L^{p}_{cyl}(\mathbb{R}^{n})=\overline{C_{0,cyl}^{\infty}(\mathbb{R}^{n})}^{L^{p}(\mathbb{R}^{n})} , \quad H^{s}_{cyl}(\mathbb{R}^{n})= \overline{C_{0,cyl}^{\infty}(\mathbb{R}^{n})}^{H^{s}(\mathbb{R}^{n})}.$$
    Here $C^{\infty}_{0,cyl}(\mathbb{R}^{n})$ is the space of compactly supported smooth functions with cylindrical symmetry. We denote by $\mathcal{H}^\beta_{cyl}(f)$ the closed subspace of all cylindrically symmetric functions in $\mathcal{H}^\beta(f)$.
\item For the case of compact manifolds with symmetry, we refer to the action of a closed subgroup $G \leq \text{Isom}_{g}(M)$ on $(M,g)$,  where $\text{Isom}_{g}(M)$ is the group of all Riemannian isometries of $M$. We work with the fractional Sobolev space of functions invariant under the action of $G$, formally defined in \cite{Hebey 2} by
    \begin{equation}
        H^{s,p}_{G}(M)=\{u \in H^{s,p}(M) \mid \forall \sigma \in G, u \circ \sigma = u\}, \quad s \in \mathbb{R}_{>0},
    \end{equation}
    where $H^{s,p}(M)$ is defined in \cite[p.689]{Banica} as 
    \begin{equation}\label{espacios de sobolev en version laplaciano R1 para intro}
        H^{s,p}(M)=\{u \in L^{p}(M) \mid \exists h \in L^{p}(M), u=(I+\Delta_{g})^{-s/2}h \}.
    \end{equation}
Furthermore, denoting by $0 \leq \lambda_1 < \lambda_2 \leq \dots \to \infty$ the discrete spectrum of $\Delta_g$ and by $\{\omega_i\}_{i=1}^{\infty}$ an associated orthonormal basis of eigenfunctions in $L^2(M,g)$, we define the space
\begin{equation}
\widetilde{H}_{f}(M,g) = \left\{h=\sum_{i=1}^{\infty} \alpha_{i}\omega_{i} ~\Big|~ \sum_{i=1}^{\infty}|\alpha_{i}|^{2}(1 + f(\lambda_{i}))^{2} < \infty \right\}.
\end{equation}
 We assume that $f$ satisfies appropriate growth conditions, and introduce the space  $\widetilde{H}_{f,G}(M,g)$, which will allow us to work with the operator $f(\Delta_{g}) + I$. This space is defined as
\begin{equation}\label{espacio de Sobolev en verson sumas con acción para intro}
    \widetilde{H}_{f, G}(M,g) = \left\{ h \in \widetilde{H}_{f}(M,g) \mid h \circ \sigma=h, \forall \sigma \in G \right\}.
\end{equation}
Similarly we consider $L^{p}_{G}(M)$ defined by $L^{p}_{G}(M)=\left\{u \in L^{p}(M) \mid u \circ \sigma=u, \forall \sigma \in G \right\}$.
\end{enumerate}

\justify

We summarize our main theorems below. First, for cylindrical symmetry on $\mathbb{R}^n$, motivated by the 
main tool for proving a variant of Theorem \cite[t.3.5]{Gorka 6}, we use the Lion's compactness theorem proved in \cite[t.III.3]{Lions}. We obtain the following 

\begin{theorem}[\ref{teorema para Rn con simetria cilindrica}]
    Suppose that $0 < \gamma < 1$ is a constant\,\footnote{We can consider the condition $\alpha > 1$ as \cite{Gorka 6}, instead of the condition $0 < \gamma < 1$, because $\gamma = \frac{\alpha - 1}{\alpha}$.
}, $f \in \mathcal{G}^\beta$ for some $\beta > \frac{n\gamma}{2}$, and that $U:\mathbb{R}^{n} \times  \mathbb{R} 
\to \mathbb{R}$ is cylindrically symmetric in $x$, is $C^{1}$ in the second variable and $U(\cdot, 0) \not\equiv 0$.  
Assume moreover that there exist functions $h \in L^2(\mathbb{R}^n)$, $g \in L^{2/\gamma}(\mathbb{R}^n)$ and some constant $C>0$ such that the following two inequalities hold,
\begin{equation}\label{desigualdades caso cilindrico para intro}
 |U(x, y) + y| \leq C(|h(x)| + |y|^{\frac{1}{1-\gamma} } ),
\end{equation}
\vspace{-0.3cm}
\begin{equation}\label{desigualdades caso cilindrico 2 para intro}
    \left|\frac{\partial}{\partial y} (U(x, y) + y)\right| \leq C(|g(x)| + |y|^{\frac{\gamma}{1 - \gamma} }).
\end{equation}
Then there exist $0 < \tau < 1$ and $0 < \rho_{\tau } < 1$ such that, whenever $\|h\|_{L^2(\mathbb{R}^n)} < \rho_{\tau}$, there exists a solution $u \in \mathcal{H}^\beta_{cyl}(f)$ to the equation $f(\Delta)u - U(\cdot, u) = 0$ with $\|u\|_{L^{2/(1-\gamma) }_{cyl}(\mathbb{R}^n)} \leq \tau $.
\end{theorem}

The most technical result of this work is an existence theorem for smooth solutions to nonlocal equations on 
compact Riemannian manifolds under the action of a closed subgroup 
$G \leq \mathrm{Isom}_g(M)$. The main tool is the following new compact 
embedding theorem, which extends the integer-order result of Hebey-Vaugon 
\cite[t.9.1]{Hebey 2} to the fractional Sobolev spaces $H^{s,q}_G(M)$ of 
$G$-invariant functions.

\begin{theorem}[\ref{Lema principal caso con accion}]
 Let $(M,g)$ be a compact, connected, oriented Riemannian $n$-manifold without boundary. Let $q>1$, $s>0$ be real numbers, and let $G$ be a closed subgroup of $\mathrm{Isom}_g(M)$. Assume that for any $x \in M$, $|\mathcal{O}_{G}(x)|=+\infty$ and set $k=\min_{x\in M}\dim \mathcal{O}_{G}(x).$ Then, $k\ge 1$ and for every $p \in [1,\infty)$ satisfying
$
s>(n-k)\left(\frac{1}{q}-\frac{1}{p}\right),$
the embedding
\begin{equation}
   H^{s,q}_{G}(M)\hookrightarrow L^{p}(M) 
\end{equation}
is compact. Equivalently, if $qs \geq n-k$, the embedding is compact for every real $p \in [1, 
\infty)$. More precisely, if $qs > n-k$, $H^{s,q}_{G}(M)$ embeds compactly into $C^{0}(M)$(and hence into $L^{\infty}(M)$).
If $qs<n-k$, the embedding is compact for $1 \leq p < \frac{(n-k)q}{n-k-sq}$, and continuous for $1 \leq p \leq \frac{(n-k)q}{n-k-sq}$. 
\end{theorem}

Using Theorem~\ref{Lema principal caso con accion} in place of the classical Rellich--Kondrachov theorem, together with Schauder's fixed-point theorem, we prove the following main existence result.

\begin{theorem}[\ref{Teorema Final}]
   Let $(M,g)$ be a compact, connected and oriented Riemannian $n$-manifold without boundary. Let $G$ be a closed subgroup of $\mathrm{Isom}_{g}(M)$ whose action on $M$ is non-transitive. 
    Suppose that $|\mathcal{O}_{G}(x)|=+\infty, \forall x \in M$, and $k=\min_{x \in M}\mathrm{dim}~\mathcal{O}_{G}(x)$. Assume that $f$ satisfies (\ref{codicion de elipticidad de una fucnion real valuada para Riemannian manifolds}),
    $U: M \times \mathbb{R} \to \mathbb{R}$ is a Carathéodory non-linearity which is $G$-invariant with $U(\cdot, 0) \not\equiv 0$ and
    $C^{1}$ in the second variable, 
    and there exist constants $\alpha>1$, $\beta \geq 1$, $C > 0$, and a function $h \in L_{G}^2(M,g)$ such that the following holds:
\begin{equation}\label{desigualdad de crecimiento intro}
|U(x, y) + y| \leq C \left( |h(x)| + |y|^\alpha \right),\quad \left| \frac{\partial}{\partial y} \left( U(x, y) + y \right) \right| \leq C ( 1 + |y|^\beta ).
\end{equation}
Then, there exists $\varepsilon \in (0,1]$ and $0<\rho_{\varepsilon} < 1$ such that whenever $\norm{h}_{L^2(M)} < \rho_{\varepsilon}$, there exists a solution $\phi \in \widetilde{H}_{f,G}(M,g)$ to the non-linear problem \begin{equation}
f(\Delta_g)\phi - U(x, \phi) = 0.
\end{equation}
\end{theorem}

The article is structured as follows: in Section \hyperref[Preliminaries and main results]{(II)}, we present the necessary preliminaries regarding nonlocal theory, both in Euclidean spaces and on compact Riemannian manifolds. We discuss the essential functional-analytic tools, review foundational previous results from \cite{Gorka 6}, and establish the rigorous representation of the operators $f(\Delta)$ and $f(\Delta_{g})$ via functional calculus and spectral theory, which subsequently allows us to define the appropriate functional spaces for our problem. In Section \hyperref[Nonlocal equations on Bounded Domains]{(III)}, we focus on the Euclidean setting and prove Theorem \ref{teorema para Rn con simetria cilindrica}, establishing the existence of solutions under cylindrical symmetry.

Finally, in Section \hyperref[seccion de compact Riemannian manifolds]{(IV)}, we 
study nonlocal equations on compact Riemannian manifolds in the presence of symmetry. First we prove the generalized fractional compact embedding Theorem \ref{Lema principal caso con accion} for $G$-invariant functions, and subsequently utilize this crucial compactness result alongside Schauder's fixed-point theorem \cite[c.11.2]{Trudinger} to establish the main existence result, Theorem \ref{Teorema Final}, for compact Riemannian manifolds with symmetries.

\newpage

\section{Preliminaries}\label{Preliminaries and main results}

We introduce appropriate function spaces to interpret the nonlocal term that appears in Equation (\ref{ecuacion no local}). Specifically, our goal in this section is to give meaning to the formal operator
\begin{equation}\label{operador no local}
    L = f(\Delta) + \text{Id}
\end{equation}
where $f$ is a measurable function, as an operator defined on (a subspace of) 
$L^2(\mathbb{R}^n)$. We will call $f$ the "symbol"; however, it need not belong to Hörmander's class of symbols\cite{Hormander}, since the framework developed here is more general and self-contained. Our approach is based on a direct functional calculus via the Fourier transform, which does not require $f$ to be smooth nor to satisfy the polynomial growth estimates characteristic of Hörmander's theory. Instead, it will suffice that $f$ satisfies the following conditions, \\
$(P)$~~\text{The function} $s \mapsto f(-s^{2})$~\text{is nonnegative}.\\
$(E_{\beta})$ There exist real numbers $\beta, R, M > 0$ such that
\begin{equation*}\label{elliptic condition over f}
M(1+|\xi|^{2})^{\frac{\beta}{2}} \leq f(-|\xi|^{2})~~,\forall \xi \in \mb{R}^{n}~\text{such that}~ |\xi|> R.
\end{equation*}
Functions satisfying these conditions are called $\mathcal{G}^{\beta}$-symbols. The condition $(E_{\beta})$ is the ellipticity condition for a pseudo-differential operator.

\begin{definition}
    Fix $f \in \mathcal{G}^{\beta}$. The space $\mathcal{H}^{\beta}(f)$ is given by all real-valued functions $g$ on $\mathbb{R}^n$ such that $g$ is measurable, the Fourier transform $\mathcal{F}(g)$ exists, and
    \begin{equation}
     \int_{\mathbb{R}^n} \left[ 1 + f\left(-|\xi|^2\right) \right]^2 \left| \mathcal{F}(g)(\xi) \right|^2 \, d\xi < \infty.
    \end{equation}
\end{definition}

\justify
This space satisfies $\mathcal{H}^{\beta}(f) \subseteq L^2(\mathbb{R}^n).$ The reason for defining this space is explained by the following computation. If we assume that $f$ is a real-analytic function, one can consider $f$ as a power series centered at $x=0$, which we may evaluate at $\Delta$, giving the formal expression
\begin{equation*}
f(\Delta)u=\sum_{n=0}^{\infty}\frac{f^{(n)}(0)}{n!}\Delta^{n}u 
\end{equation*}
Applying the Fourier transform to this equation yields
\begin{align*} 
\mathcal{F}(f(\Delta)u)(\xi)&=\sum_{n=0}^{\infty}\frac{f^{(n)}(0)}{n!}\mathcal{F}(\Delta^{n}u)= \sum_{n=0}^{\infty}\frac{f^{(n)}(0)}{n!}(-|\xi|^{2})^{n}\mathcal{F}(u) = f(-|\xi|^{2})\mathcal{F}(u)(\xi) 
\end{align*}

Therefore we define the operator $f(\Delta)u$ by
\begin{equation}\label{forma de operador f}
    f(\Delta)u = \mathcal{F}^{-1}(f(-|\xi|^2)\mathcal{F}(u)(\xi)).
\end{equation}
If $u \in \mathcal{H}^{\beta}(f)$, then the function $f(-|\xi|^2)\mathcal{F}(u) \in L^2(\mathbb{R}^n)$, and $f(\Delta)u$ is well defined, so that $f(\Delta)$ is a linear operator from $\mathcal{H}^{\beta}(f)$ to $L^2(\mathbb{R}^n)$.

\begin{definition}\label{espacio Hb permite la interpretacion de L}
    Let $f$ be a $\mathcal{G}^{\beta}$-symbol. The operator $L = f(\Delta) + \text{I}$ acts on functions $u \in \mathcal{H}^{\beta}(f)$ as follows
    \begin{equation}\label{operador no local L}
        Lu = \mathcal{F}^{-1} \left( \mathcal{F}(u) + f(-|\xi|^2)\mathcal{F}(u) \right)
    \end{equation}
Since $u \in \mathcal{H}^{\beta}(f)$, we have $Lu \in L^2(\mathbb{R}^n)$.
\end{definition}

A function used repeatedly is
\begin{equation}\label{funcion K}
 \widehat{K}(\xi)=\cfrac{1}{1+f(-|\xi|^{2})}   
\end{equation}
In \cite[l.2.1]{Gorka 6} it is proved that for $f$ a $\mathcal{G}^{\beta}$-symbol with $\beta > n/2$, one has $\widehat{K} \in L^{2}(\mathbb{R}^{n})$. This implies there exists a function $K_{f}$ such that $\mathcal{F}(K_{f})=\widehat{K}.$

\begin{proposition}\label{espacio Hillbert Hb}
A measurable function $g$ belongs to the space $\mathcal{H}^{\beta}(f)$, with $\beta > n/2$, if and only if $g = K_f * \tilde{g}$ for some $\tilde{g} \in L^2(\mathbb{R}^n)$. Moreover, if we endow $\mathcal{H}^{\beta}(f)$ with the inner product
\begin{equation}
    \langle g_1, g_2 \rangle = \int_{\mathbb{R}^n} \left( 1 + f(-|\xi|^2) \right)^2 \mathcal{F}(g_1)(\xi) \overline{\mathcal{F}(g_2)(\xi)} \, d\xi
,
\end{equation}
then $\mathcal{H}^{\beta}(f)$ is a Hilbert space and the map $\mathcal{K} : L^2(\mathbb{R}^n) \to \mathcal{H}^{\beta}(f)$ given by $ \mathcal{K}(g) = K_f * g$\label{isomorfismo isometrico entre los espacios Hb}
, is an isometric isomorphism (a bijective linear isometry).
\end{proposition}

\begin{lemma}\label{Lemma de contencion epscios Hb}
   \cite[l.2.2]{Gorka 6}. Let $f$ be a $\mathcal{G}^{\beta}$-symbol. We have:\\
i) For each $s \in \mathbb{R}$ with $s \leq \beta$, the embedding 
\begin{equation}
  \mathcal{H}^{\beta}(f) \hookrightarrow H^{s}(\mathbb{R}^n)
\end{equation}
holds. In other words,
$
\| u \|_{H^{s}(\mathbb{R}^n)} \leq C(s, \beta) \| u \|_{\mathcal{H}^{\beta}(f)},
$
for some constant $C(s, \beta) > 0$.\\
ii) For each integer $k \geq 1$ such that $\frac{n}{2} + k < s \leq \beta$, one has $$\mathcal{H}^{\beta}(f) \hookrightarrow C^{k}(\mathbb{R}^n).$$
iii) For each $0 < \alpha < 1$ such that $\frac{n}{2} + \alpha = s \leq \beta$, the embedding $$\mathcal{H}^{\beta}(f) \hookrightarrow C^{\alpha}(\mathbb{R}^n)$$
holds, where $C^{\alpha}(\mathbb{R}^n)$ is the space of H\"{o}lder continuous functions of order $\alpha$.
\end{lemma}

\vspace{0.5cm}

\subsection{The Linear Equation}

\begin{definition}
    A function $u : \mathbb{R}^n \to \mathbb{R}$ is \emph{real-analytic} if it is infinitely differentiable and can be represented locally by a convergent power series. Thus, for each point $z \in \mathbb{R}^n$ there exists an open neighborhood $U$ of $z$ such that
\begin{equation}
    u(x) = \sum_{\alpha} \frac{ \partial^\alpha u(z)}{\alpha!} (x - z)^\alpha,
\end{equation}
the sum runs over all multi-indices. The space of real-analytic functions is denoted $C^\omega(\mathbb{R}^n)$. 
\end{definition}

 A function $u(z)$ is real-analytic if the \emph{remainder term} tends to zero. In fact $|R_{m}(x)| \to 0$, where
\begin{equation*}\label{remainder term 2}
    R_{m}(x)=\sum_{|k|=m}\frac{\partial^{k}u(z+t(x-z))(x-z)^{k}}{k!}, \quad t \in [0,1].
\end{equation*}

\begin{theorem}\label{teorema caso lineal}\cite[t.3.1]{Gorka 6}. 
    Let $f$ be a $\mathcal{G}^{\beta}$-symbol. For each $g \in L^2(\mathbb{R}^n)$ there exists a unique solution $u_g \in \mathcal{H}^\beta(f)$ to the linear problem
\begin{equation}\label{ecuacion lineal nolocal}
    L u = g.
\end{equation}
Moreover, the following identity holds
\begin{equation}
  \|u_g\|_{\mathcal{H}^\beta(f)} = \|g\|_{L^2(\mathbb{R}^n)}.  
\end{equation}
\end{theorem}

 In the case $\beta > n/2$, the solution $u_{g}$ is given by the expression
\begin{equation}\label{forma lineal de la solución}
    u_g = \mathcal{F}^{-1} \left( \frac{\mathcal{F}(g)}{1 + f(-|\xi|^2)} \right) = K_{f} \ast g.
\end{equation}

\begin{theorem}\label{funcion real-analitica}
    Suppose $f \in \bigcap_{\beta > 0}\mathcal{G}^{\beta}$ and that there exists $c>0$ such that
    \begin{equation}
        \lim_{m \to \infty} \frac{c^{m}}{m!}\int_{0}^{\infty}\frac{r^{2m+n-1}}{[1+f(-r^{2})]^{2}}dr=0.
    \end{equation}
Then the solution $u$ of the linear problem $L u = g$ with $g \in L^2(\mathbb{R}^n)$ satisfies $u \in C^\omega(\mathbb{R}^n)$.
\end{theorem}

\newpage

\subsection{The Nonlinear Equation}

\justify

In this section we study the nonlinear equation
\begin{equation}\label{ecuacion basica no-lineal}
    f(\Delta)u - U(\cdot , u)=0
\end{equation}
Before analyzing this equation in the general case, let us consider a nonlocal term of the form $U(x,y)=y + \varphi(x)V(x,y)$, where $x \in \mathbb{R}^{n}$ and $y \in \mathbb{R}$, and $\varphi \in C^{\infty}_{0}(\mathbb{R}^{n})$. This kind of equation is called \emph{localized case},

\begin{theorem}\label{caso localizado}\cite[t.3.3]{Gorka 6}.
    Suppose $f \in \mathcal{G}^\beta$ with $\beta > \frac{n}{2}$, $n > 1$, and consider the function $U_\delta$, with $\delta > 0$, given by:
\begin{equation}\label{primera ecuacion no lineal}
U_\delta(x, y) = -y + \delta \varphi(x)V(x, y),
\end{equation}
where $\varphi \in C_0^\infty(\mathbb{R}^n)$ and $V$ is a $C^1$-function on $\mathbb{R}^n \times \mathbb{R}$. Assume there exist a constant $\alpha \geq 1$, a constant $C > 0$ and a function $h \in L^2(\mathbb{R}^n)$ such that the following inequalities hold:
\begin{equation}\label{primera desigualdad caso localizado}
|V(x, y)| + \left|\frac{\partial}{\partial x_i} V(x, y)\right| \leq C(h(x) + |y|^\alpha), \quad i = 1, \ldots, n,
\end{equation}
\begin{equation}\label{segunda desigualdad caso localizado}
\left|\frac{\partial}{\partial y} V(x, y)\right| \leq C(1 + |y|^\alpha).
\end{equation}
Then, for $\delta$ sufficiently small, there exists a solution $u \in \mathcal{H}^\beta(f)$ of the problem (\ref{ecuacion basica no-lineal}) with the nonlinearity $U$ given by (\ref{primera ecuacion no lineal}).
\end{theorem}

\begin{corollary}
    Suppose $U$ is given by (\ref{primera ecuacion no lineal}), the hypotheses of Theorem \ref{funcion real-analitica} are satisfied, and the growth assumptions of Theorem \ref{caso localizado} hold. Then there exists a $C^\omega$ solution to the problem (\ref{ecuacion basica no-lineal}).
\end{corollary}

Now we work with a nonlinearity $U : \mathbb{R}^{n} \times  \mathbb{R} \to \mathbb{C}$ which is \emph{spherically symmetric} in $x$, which means that the nonlinearity satisfies $U(|x|,y)=U(x,y)$.

\begin{definition}
\begin{equation}
    L^2_r(\mathbb{R}^n) = \overline{C_{0,r}^\infty(\mathbb{R}^n) }^{L^2(\mathbb{R}^n)}~~,~~
H^s_r(\mathbb{R}^n) = \overline{C_{0,r}^\infty(\mathbb{R}^n) }^{H^s(\mathbb{R}^n)},
\end{equation}  
where $C_{0,r}^\infty(\mathbb{R}^n)$ denotes the space of smooth radial compactly supported functions on $\mathbb{R}^n$. Denote by $\mathcal{H}^\beta_r(f)$ the closed subspace of radial functions in $\mathcal{H}^\beta(f)$.
\end{definition}

\begin{theorem}\label{problema lineal para el caso cilindrico}
    If $g \in L^2_r(\mathbb{R}^n)$, then the unique solution $u \in \mathcal{H}^\beta(f)$ of the problem \eqref{ecuacion lineal nolocal} is a radial function.  
Moreover, the identity holds:
\begin{equation}
\|u\|_{\mathcal{H}^\beta_r(f)} = \|g\|_{L^2_r(\mathbb{R}^n)}.
\end{equation}
\end{theorem}

\begin{theorem}\label{Caso Radial}\cite[t.3.5]{Gorka 6}.
    Suppose $\alpha > 1$ is a constant, $f \in \mathcal{G}^\beta$ for some $\beta > \frac{n}{2} \left(\frac{\alpha - 1}{\alpha}\right)$, and $U$ is spherically symmetric in $x$.  
Assume moreover that there exist functions $h \in L^2(\mathbb{R}^n)$ and $g \in L^{\frac{2\alpha}{\alpha -1} }(\mathbb{R}^n)$ such that the following two inequalities hold:
\begin{equation}\label{primera desigualdad caso radial}
|U(x, y) + y| \leq C(|h(x)| + |y|^\alpha),
\end{equation}
\begin{equation}\label{segunda desigualdad caso radial}
    \left|\frac{\partial}{\partial y} (U(x, y) + y)\right| \leq C(|g(x)| + |y|^{\alpha - 1})
\end{equation}
for some constant $C > 0$. Then there exist $0 < \varepsilon < 1$ and $0 < \rho_{\varepsilon} < 1$ such that, whenever $\|h\|_{L^2(\mathbb{R}^n)} < \rho_{\varepsilon}$, there exists a radial solution $u \in \mathcal{H}^\beta_r(f)$ to the equation
\begin{equation}
f(\Delta)u - U(\cdot, u) = 0
\end{equation}
with $\|u\|_{L^{2\alpha}_r(\mathbb{R}^n)} \leq \varepsilon$.
\end{theorem}

\begin{corollary}
    Suppose $U$ is spherically symmetric in $x$, the hypotheses of Theorem \ref{funcion real-analitica} are satisfied, and the growth assumptions of Theorem \ref{Caso Radial} hold. Then there exists a $C^\omega$ solution to the problem \ref{ecuacion basica no-lineal}.
\end{corollary}

\vspace{0.12cm}

\subsection{Nonlocal Equations on Compact Riemannian Manifolds}\label{Nonlocal Equations on Compact Riemannian Manifolds}

\justify
In this section we consider nonlocal equations of the form  
\begin{equation}\label{ecuacion no local en Riemannian manifolds}
  f(\Delta_g)\phi - U(x, \phi) = 0  
\end{equation}
on a compact, connected, oriented Riemannian manifold $(M,g)$. The \emph{Laplacian} $\Delta_g$ associated to $g$ is a second-order differential operator acting on smooth functions, i.e., $\Delta_g: C^\infty(M) \to C^\infty(M)$, which is given by
\begin{equation}
 \Delta_g f =- \text{div}_g(\nabla_g f),   
\end{equation}
where $\nabla_g f$ is the gradient of $f$, defined as the unique vector field $X$ such that 
$
g(X, Y) = df(Y),~\forall Y \in \mathfrak{X}(M)
$. The operator $\text{div}_g$ is the divergence of a vector field $X=X^{i}\frac{\partial}{\partial x^{i}}$ in a local coordinates $(U,(x^{i}))$, defined in terms of the volume form $dV_g$ induced by $g$,
\begin{equation}
    \text{div}_g(X) = \frac{1}{\sqrt{|g|}} \partial_i \big(\sqrt{|g|} X^i\big),
\end{equation}
where $\det g$ is the determinant of the matrix $(g_{ij})$ of the metric in local coordinates. 
In local coordinates $(x^i)$ on $M$, the Laplacian $\Delta_g$ of a function $f$ is written as,
\begin{equation}\label{Operador de Laplace-Beltrami}
    \Delta_g f =  -\frac{1}{\sqrt{|g|}} \partial_i \left(\sqrt{|g|} g^{ij} \partial_j f \right),
\end{equation}
where $g^{ij}$ are the coefficients of the inverse of the metric $g_{ij}$.
Another intrinsic definition of $\Delta_g$ is
$
\Delta_{g} u = -\ast d \ast d u,  
$
where $d$ denotes the exterior derivative and $\ast$ is the Hodge star operator. Recall that by definition the Hodge star operator $\star: \Omega^k(M) \to \Omega^{n-k}(M)$ assigns to any $k$-form $\eta \in \Omega^k(M)$ the $(n - k)$-form $\star\eta \in \Omega^{n-k}(M)$ such that 
\begin{equation*}
   \omega \wedge \star\eta = \langle \omega, \eta \rangle \omega_g, \quad  \omega \in \Omega^k(M).
\end{equation*}

\justify
We recall the following definition of spaces $(L^{2}(M,g), dV_{g})$

\begin{definition}
    The space $L^2(M, g)$ is defined as the set of square-integrable functions on the Riemannian manifold $(M, g)$ with respect to the measure induced by the metric $g$, namely
\begin{equation}\label{espacio de funciones cuadrado integrables en una variedad Riemanniana}
    L^2(M, g) = \left\{ f: M \to \mathbb{R}(\mathbb{C}) ~\Big|~ \int_M |f|^2 \, dV_g < \infty \right\}
\end{equation}
Here $dV_g$ is the volume measure on $M$ associated with $g$, defined by $dV_g = \sqrt{\det(g)} \, dx^1 \wedge dx^2 \wedge \dots \wedge dx^n$ (also denoted $\omega_{g}$). This is a Hilbert space equipped with the inner product
$$
\langle f, h \rangle_{L^2(M)} = \int_M f \overline{h} \, dV_g,
$$
and the norm induced by this inner product is
  $$
    \|f\|_{L^2(M)} = \sqrt{\langle f, f \rangle_{L^2(M)}} = \left( \int_M |f|^2 \, dV_g \right)^{1/2}.
    $$
\end{definition}

\justify

The Laplace–Beltrami operator $\Delta_{g}$ is originally defined on $C^{\infty}(M)$, but it extends naturally to $L^{2}(M,g)$ from a dense domain; indeed
\begin{equation*}
\Delta_{g}: \mathcal{D}(\Delta_{g}) \subset L^{2}(M,g) \to L^{2}(M,g)
\end{equation*}
where the maximal domain $\mathcal{D}(\Delta_{g})=\{u \in H^{2}_{\text{loc}}(M) \cap L^{2}(M) \mid \Delta_{g}u \in L^{2}(M) \}$.

\begin{proposition}\cite{Backman,Canzani,Rosenberg}.
    The operator $\Delta_{g}$ has the following properties.\\
    i)~ $\Delta_g$ is a second-order elliptic operator and is essentially self-adjoint with respect to the inner product on $L^2(M, g)$, that is
$
\langle f, \Delta_g h \rangle = \langle \Delta_g f, h \rangle.
$\\
ii)~ If $M$ is compact, the spectrum of $\Delta_g$ is a discrete set of nonnegative eigenvalues $0\leq \lambda_{1}<\lambda_{2} \leq \lambda_{3} \leq \ldots$, with an orthonormal basis of eigenfunctions $\{\omega_i\}$ in $L^2(M, g)$ such that $
\Delta_g \omega_i = \lambda_i \omega_i.$\\
iii)~ The operator $\Delta_g$ commutes with isometries of the manifold $(M, g)$. This means that if $\varphi: M \to M$ is an isometry then
$$
\Delta_g (f \circ \varphi) = (\Delta_g f) \circ \varphi.
$$
iv)~ The operator $\Delta_{g}$ is unbounded and $\Delta_{g}\geq 0$.
\end{proposition}

\justify

The operator $\Delta_{g}: H^{2}(M,g) \to L^{2}(M,g)$ is an unbounded self-adjoint operator on a compact Riemannian manifold, and the resolvent operator $R_{\lambda}(\Delta_{g})$ defined by
\begin{equation*}
R_{\lambda}(\Delta_{g}): L^{2}(M,g) \to H^{2}(M,g), \quad \psi \mapsto  (\Delta_{g} - \lambda I)^{-1}\psi
\end{equation*}
is linear, bounded, self-adjoint and compact \cite{ReedI}.
Therefore, by the Theorem \cite[p.5.12]{Schmudgen} 
there exist $(\lambda_{n})_{n \in \mathbb{N}} \subset \mathbb{R}$ and an orthonormal basis $\{ \omega_{n} \mid n \in \mathbb{N}\}$ of $L^{2}(M,g)$ such that $0\leq\lambda_{1}<\lambda_{2} \leq \lambda_{3} \leq \ldots$ with finite multiplicity and 
accumulation at infinity, i.e., $|\lambda_{j}| \to \infty$, and $\Delta_{g}\omega_{n}=\lambda_{n} \omega_{n}$ for $n \in \mathbb{N}$. 
Another way to justify the existence of an orthonormal basis of $L^{2}(M,g)$ consisting of eigenfunctions of $\Delta_{g}$
is via the Hodge Theorem for functions \cite[t.1.29]{Rosenberg}. 
Using the spectral theorem for unbounded self-adjoint operators 
\cite[t.5.7]{Schmudgen} there exists a
unique spectral measure $E_{\Delta_{g}}: \mathcal{B}(\mathbb{R}) \to \mathcal{L}(L^{2}(M,g))$ such that
\begin{equation}\label{medida espectral para el laplaciano beltrami}
\hspace{0.3cm}~\Delta_{g}=\int_{\mathbb{R}}\lambda dE_{\Delta_{g}}(\lambda).
\end{equation}
In order to use the functional calculus, by the assignment $f \mapsto f(\Delta_{g})$, one requires that $f \in \mathcal{S}(\mathbb{R}, \mathcal{B}(\mathbb{R}), E_{\Delta_{g}})$, which is the set of all $E_{\Delta_{g}}$-a.e. finite Borel functions $f: \mathbb{R} \to \mathbb{C} \cup \{\infty\}$, that is, $E_{\Delta_{g}}(\{t \in \mathbb{R} \mid f(t)=\infty \})=0$. Since $f \in \mc{G}^{\beta}$, and $f$ is a real-valued measurable function, it follows that $f \in \mathcal{S}(\mathbb{R}, \mathcal{B}(\mathbb{R}), E_{\Delta_{g}})$. Applying that $\sigma(\Delta_{g})$ is purely discrete, we have that
\begin{equation}\label{expresion 2, Gorka, On a class}
    f(\Delta_{g})=f\left( \int_{\mathbb{R}}\lambda dE_{\Delta_{g}}(\lambda)    \right)= \int_{\mathbb{R}}f(\lambda)dE_{\Delta_{g}}(\lambda)=\sum_{i=1}^{\infty}f(\lambda_{i})
E_{\Delta_{g}}(\{\lambda_i\})\footnote{Recall that $E_{\Delta_{g}}(\lambda_{i})=E_{\Delta_{g}}((-\infty, \lambda_{i}])$ is a resolution of the identity.}
\end{equation}
For an operator with a point spectrum, the spectral measure $E_{\Delta_{g}}(\cdot)$ evaluated at the singleton $\{\lambda_{i}\}$ yields $P_{i}$. Specifically, each $P_{j}$ is the orthogonal projection operator onto the closed subspace corresponding to the eigenspace associated with the eigenvalue $\lambda_{j}$. More concretely
\begin{equation}
E_{\Delta_{g}}(\{\lambda_{i}\}) =P_{i}, \quad E_{\Delta_{g}}((-\infty, \lambda_{i}])=\sum_{\lambda_{j} \leq \lambda_{i}}P_{j}.
\end{equation}
Since every closed subspace $\mathcal{H} \subseteq L^{2}(M,g)$ 
can be written as $\mathcal{H}=\overline{\text{span}\{\omega_{n_{i}}\}_{i \in I}}$, it 
follows that for each $C \in \mathcal{B}(\mathbb{R})$ the expression $E_{\Delta_{g}}(C)\psi$ has the form
$$E_{\Delta_{g}}(C)\psi=\sum_{\lambda_{i} \in C}\langle \psi, \omega_{n_{i}}   \rangle \omega_{n_{i}}$$
and therefore (\ref{expresion 2, Gorka, On a class}) is equivalent to
\begin{equation}
 f(\Delta_{g})\psi= \int_{\mathbb{R}}f(\lambda)dE_{\Delta_{g}}(\lambda)\psi = \sum_{i=1}^{\infty}f(\lambda_{i})P_{i}\psi   =\sum_{i=1}^{\infty}f(\lambda_{i})\langle \psi, \omega_{i}  \rangle \omega_{i}  = \sum_{i=1}^{\infty}f(\lambda_{i})\psi_{i}\omega_{i}.
\end{equation}

\justify

The last identity will serve as a tool to study the linear equation
\begin{equation}\label{1.(33)}
    f(\Delta_{g})\phi + \nu \phi=h, \quad  \nu \in \mb{R}_{>0}
\end{equation}
where $f$ is a nonnegative function. If $h \in L^2(M,g)$, using the eigenfunction basis of $\Delta_{g}$ one can write $h = \sum h_{i}\omega_{i}$, and a straightforward calculation yields $\phi_i = \frac{h_i}{f(\lambda_i) + \nu}.$ Hence $\phi$ is
\begin{equation}\label{unica solucion caso lineal caso Riemanniano}
    \phi = \sum_{i=1}^\infty \frac{h_i}{f(\lambda_i) + \nu} \omega_i,
\end{equation}
By its form this function belongs to $L^{2}(M,g)$; however, it need not be smooth for that we introduce the following spaces. 
Consider the following Hilbert space $(l^{2,k}(\mathbb{N}),  \langle \alpha , \beta  \rangle_{l^{2,k}(\mathbb{N})} )$
\begin{equation*}
    l^{2,k}(\mathbb{N}) = \left\{ \alpha = (\alpha_i) ~\Big|~ \alpha_{i} \in \mathbb{R}, \sum_{i=1}^{\infty} |\alpha_i|^2 (1 + \lambda_i^2)^k < \infty \right\}, ~
    \langle \alpha , \beta  \rangle_{l^{2,k}(\mathbb{N})} = \sum_{i=1}^{\infty}\alpha_{i}\beta_{i}(1+\lambda_{i}^{2})^{k}
\end{equation*}
The notation with the eigenvalues $\lambda_{i}$ squared is used in \cite{Gorka 6} to emphasize that the eigenvalues are positive. Also define the space
\begin{align}
\widetilde{H}^k(M,g) &=  \left\{ h=\sum_{i=1}^{\infty}\alpha_{i}\omega_{i} ~\Big|~ \sum_{i=1}^{\infty}|\alpha_{i}|^{2}(1+\lambda_{i}^{2})^{k} < \infty    \right\}
\end{align}
This space is endowed with the following inner product
\begin{equation}
    \langle  h_{1} , h_{2}   \rangle_{\widetilde{H}^k(M,g)} = \langle (\alpha_{i})  , (\beta_{i})  \rangle_{l^{2,k}(\mathbb{N})} = \sum_{i=1}^{\infty}\alpha_{i}\beta_{i}(1+\lambda_{i}^{2})^{k},
\end{equation}

\justify

The space $\widetilde{H}^{k}(M,g)$ is isometrically isomorphic to $l^{2,k}(\mathbb{N})$. Moreover, if $k \in \mathbb{N}$ then $\widetilde{H}^{k}(M,g) = H^{2k}(M,g)$ \cite[l.4.1]{Gorka 6}, where the latter is the Sobolev space on the Riemannian manifold $(M,g)$.\\
We assume that $f$ is a nonnegative real-valued function satisfying the ellipticity condition for the spectrum $(\lambda_{i})$ of the Laplace operator, namely
\begin{equation}\label{codicion de elipticidad de una fucnion real valuada para Riemannian manifolds}
   \forall k>0, \exists\, C_{k}>0 :  C_{k}(1+\lambda_{i}^{2})^{k/2} \leq 1+ f(\lambda_{i}).
\end{equation}
This condition replaces the $(P)$ and $(E_{\beta})$ conditions. Introduce the space
\begin{equation}\label{espacio de Sobolev en verson sumas}
    \widetilde{H}_{f}(M,g) = \left\{\phi=\sum_{i=1}^{\infty} \alpha_{i}\omega_{i} ~\Big|~ \sum_{i=1}^{\infty}|\alpha_{i}|^{2}(1 + f(\lambda_{i}))^{2} <  \infty   \right\}
\end{equation}
Thanks to the ellipticity condition, for every $k$ we have the embedding
\begin{equation}\label{embedding caso riemaniano}
  \widetilde{H}_{f}(M,g) \hookrightarrow H^{2k}(M,g)  
\end{equation}
If $h \in  \widetilde{H}_{f}(M,g)$ then $h$ is a smooth function on $M$ by \cite[t.2.7]{Hebey 2}. 
\vspace{0.12cm}

Define the operator $L_{c}=f(\Delta_{g}) + \mathrm{Id}$ by
\begin{equation}\label{operador caso Riemanniano}
 L_{c}: \widetilde{H}_{f}(M,g) \longrightarrow L^{2}(M,g), \quad \phi \longmapsto \sum (f(\lambda_{i}) + 1)\phi_{i}\omega_{i} 
\end{equation}

\begin{lemma}\label{sistema lineal caso riemanniano}
    For every $h \in L^2(M,g)$, the linear problem
\begin{equation}
f(\Delta_g)\phi + \phi = h
\end{equation}
admits a unique solution $\phi \in \widetilde{H}_f(M,g)$. Moreover,
$\norm{\phi}_{\widetilde{H}_{f}(M,g)}=\norm{h}_{L^{2}(M,g)}$.
\end{lemma}

\vspace{0.12cm}

\begin{theorem}\label{teorema ecuacion no local caso Riemannian manifold}\cite[t.4.1]{Gorka 6}.
      Suppose $f$ satisfies (\ref{codicion de elipticidad de una fucnion real valuada para Riemannian manifolds}), that $U: M \times \mathbb{R} \to \mathbb{R}$\,\footnote{In the theorem \cite[t.4.1]{Gorka 6}, the Carathéodory hypothesis is implicit.}, and that there exist constants $\alpha > 1$, $\beta \geq 1$, $C > 0$, and a function $h \in L^2(M,g)$ such that the following inequalities hold:
\begin{equation}
|U(x, y) + y| \leq C \left( |h(x)| + |y|^\alpha \right),
\end{equation}
\begin{equation}
    \left| \frac{\partial}{\partial y} \left( U(x, y) + y \right) \right| \leq C ( 1 + |y|^\beta )
\end{equation}
Then, if $\norm{h}_{L^2(M,g)}$ is sufficiently small, there exists a solution $\phi \in \widetilde{H}_f(M,g)$ to the nonlinear problem
\begin{equation}
f(\Delta_g)\phi - U(x, \phi) = 0.
\end{equation}
\end{theorem}

\newpage

\section{Nonlinear Equations under Symmetries on \texorpdfstring{$\mathbb{R}^{n}$}{Rⁿ}}\label{Nonlocal equations on Bounded Domains}

This section extends the previous framework to encompass a broader class of symmetric solutions. While the original treatment \cite{Gorka 6} focused primarily on radial symmetry, we now consider the more general case of cylindrical symmetry, which arises naturally in numerous physical applications \cite{Calgani, Witten} including wave propagation in stratified media and pattern formation in anisotropic systems \cite{Moeller}. We investigate the nonlinear nonlocal equation
\begin{equation}\label{cylindrical eq}
f(\Delta)u - U(\cdot, u) = 0
\end{equation}
where the nonlinearity $U : \mathbb{R}^{n} \times \mathbb{R} \to \mathbb{C}$ exhibits cylindrical symmetry in the spatial variable. This represents a significant generalization of the spherically symmetric case studied previously, as cylindrical symmetry preserves invariance under orthogonal transformations within specific subspaces rather than the entire ambient space. \\
We will use the following definition for the cylindrically symmetry, that appear on Lions \cite{Lions}. 
Let $m \geq 2$ and let $n=\sum_{i=1}^{m}n_{i}$ with $n_{i} \geq 2$, where $n_{i}$ are integers. For any $(x_{1}^{0},...,x_{m}^{0}) \in \mathbb{R}^{n_{1}} \times ... \times \mathbb{R}^{n_{m}}$ and any integer $i \in \{1,...,m\}$, the function $u:\mathbb{R}^{n} \to \mathbb{R}$ is said to be \emph{cylindrically symmetric} if the function
\begin{equation}\label{función de cilindrically symmetry}
  u \circ \psi : \mathbb{R}^{n_{i}} \longrightarrow \mathbb{R}, \quad x_{i} \longmapsto u(x_{1}^{0},...,x_{i-1}^{0},x_{i},x_{i+1}^{0},..,x_{m}^{0}) 
\end{equation}
is radially symmetric, i.e.
invariant under the action of $\text{O}(n_{i})$. Here the function $\psi: \mathbb{R}^{n_{i}} \to \mathbb{R}^{n}$ take an $x_{i} \in \mathbb{R}^{n_{i}}$ and fix other coordinates $x_{j}^{0} \in \mathbb{R}^{n_{j}}$ to obtain $(x_{1}^{0},...,x_{i-1}^{0},x_{i},x_{i+1}^{0},..,x_{m}^{0})$.
Indeed, for every $A \in \text{O}(n_{i})$ we have $u \circ \psi(Ax_{i})=u \circ \psi(x_{i})$, hence
$$u(x_{1}^{0},...,x_{i-1}^{0},Ax_{i},x_{i+1}^{0},...,x_{m}^{0})=u(x_{1}^{0},...,x_{i-1}^{0},x_{i},x_{i+1}^{0},..,x_{m}^{0}).$$
This extension not only broadens the theoretical foundation of nonlocal equations but also significantly expands their applicability to physical systems where full spherical symmetry is too restrictive, while still maintaining sufficient structure to permit mathematical analysis.

\begin{definition}
    We define the spaces
    $$L^{p}_{cyl}(\mathbb{R}^{n})=\overline{C_{0,cyl}^{\infty}(\mathbb{R}^{n})}^{L^{p}(\mathbb{R}^{n})} , \quad H^{s}_{cyl}(\mathbb{R}^{n})= \overline{C_{0,cyl}^{\infty}(\mathbb{R}^{n})}^{H^{s}(\mathbb{R}^{n})}.$$
    Here $C^{\infty}_{0,cyl}(\mathbb{R}^{n})$ is the space of compactly supported smooth functions with cylindrical symmetry. Denote by $\mathcal{H}^\beta_{cyl}(f)$ the closed subspace of all cylindrically symmetric functions in $\mathcal{H}^\beta(f)$.
\end{definition}

\begin{lemma}\label{funcion cilindricamente simetrica}
    Let $f$ be a $\mathcal{G}^{\beta}$-symbol with $\beta> n/2$. Then
\begin{equation}
    \widehat{K}(\xi)=\frac{1}{1+f(-|\xi|^{2})}, \quad \xi \in \mathbb{R}^{n}
\end{equation}
    is cylindrically symmetric.
\end{lemma}

\begin{proof}
Let $m,n \in \mathbb{Z}^{+}$, $m \geq 1$ and $n  \geq 2$, such that $n=\sum_{i=1}^{m}n_{i}$, with $n_{i} \geq 2$ integers for every $i \in \{1,...,m\}$. Fix $i \in \{1,...,m\}$ and 
considered $(\xi_{1}^{0},..., \xi_{i},...,\xi_{m}^{0}) \in \mathbb{R}^{n_{1}} \times \ldots \times \mathbb{R}^{n_{i}} \times \ldots \times   \mathbb{R}^{n_{m}}$, take $A \in \text{O}(n_{i})$, and considered the computation
\begin{align*}
\widehat{K} \circ \psi(\xi_i)&=\widehat{K}(\xi_1^0,\dots,\xi_{i-1}^0,\xi_i,\xi_{i+1}^0,\dots,\xi_m^0)=\frac{1}{1+f( -| (\xi_1^0,\dots,\xi_{i-1}^0,\xi_i,\xi_{i+1}^0,\dots,\xi_m^0)|^{2})} \\
&=\frac{1}{1+f\Bigl(-\Bigl(\displaystyle\sum_{j\neq i}\|\xi_j^0\|^{2}+\|\xi_i\|^{2}\Bigr)\Bigr)} = \frac{1}{1+f\Bigl(-\Bigl(\sum_{j\neq i}\|\xi_j^0\|^{2}+\|A\xi_i\|^{2}\Bigr)\Bigr)}\\
&=\widehat{K}(\xi_1^0,\dots,\xi_{i-1}^0,A\xi_i,\xi_{i+1}^0,\dots,\xi_m^0) =\widehat{K} \circ \psi(A\xi_i).
\end{align*}
Therefore $\widehat{K}(\xi)$ is a cylindrically symmetric function.
\end{proof}

\begin{lemma} \label{teorema de cerradura}
    Let $f \in \mathcal{G}^\beta$. If $g \in L^2_{cyl}(\mathbb{R}^n)$, then the unique solution $u \in \mathcal{H}^\beta(f)$ to the linear problem $Lu=g$ is a cylindrically symmetric function.  
    Moreover, the following identity holds:
    \begin{equation}
        \|u\|_{\mathcal{H}^\beta_{cyl}(f)} = \|g\|_{L^2_{cyl}(\mathbb{R}^n)}.
    \end{equation}
\end{lemma}

\begin{proof}
We know that the solution $u$ of the linear system \eqref{teorema caso lineal} satisfies $\mc{F}(u)(\xi)=\widehat{K}(\xi)\mc{F}(g)(\xi)$. Let $m\geq 2$ and let $n=\sum_{i=1}^{m}n_{i}$, where the $n_{i}\geq 2$ are integers. Fix an arbitrary index $j \in \{1,...,m\}$ and an orthogonal matrix $R \in \text{O}(n_{j})$. Define the map $\widetilde{R}: \mb{R}^{n} \to \mb{R}^{n}$ given by
\begin{equation*}
   \widetilde{R}: \mb{R}^{n} \to \mb{R}^{n}, \quad (x_{1},...,x_{j},...,x_{m}) \mapsto (x_{1},...,Rx_{j},...,x_{m}) 
\end{equation*}
where $x_{k} \in \mb{R}^{n_{k}}$ for all $k \in \{1,\dots,m\}$; this map satisfies $\widetilde{R} \in \text{O}(n)$.
Using the property of the Fourier transform on $L^2(\mb{R}^n)$ under orthogonal transformations, namely, $\mc{F}(h \circ \mc{Q}) = \mc{F}(h) \circ \mc{Q}$ for any $\mc{Q} \in \text{O}(n)$ and $h \in L^2(\mb{R}^n)$, we have that 
$\mc{F}(g) \circ \widetilde{R}(\xi)=\mc{F}(g)(\xi)$ for $\text{a.e.}$ $\xi \in \mb{R}^{n}$, because $g$ is cylindrically symmetric.\\
Since $\widehat{K}$ is cylindrically symmetric by Lemma \ref{funcion cilindricamente simetrica}, we have that $\widehat{K}(\widetilde{R}\xi)=\widehat{K}(\xi)$ and 
$$\mc{F}(u)(\widetilde{R}\xi)=\widehat{K}(\widetilde{R}\xi)\mc{F}(g)(\widetilde{R}\xi)=\widehat{K}(\xi)\mc{F}(g)(\xi)=\mc{F}(u)(\xi)$$
for $\text{a.e.}$ $\xi \in \mb{R}^{n}$; hence $\mc{F}(u \circ \widetilde{R})=\mc{F}(u) \circ \widetilde{R}=\mc{F}(u)$. Since the Fourier transform is a unitary operator on $L^{2}(\mb{R}^{n})$, it follows that $u \circ \widetilde{R}=u$. Since the index $j$ and the orthogonal matrix $R \in \text{O}(n_j)$ were chosen arbitrarily, $u$ is cylindrically symmetric, concluding that $u \in \mathcal{H}^\beta_{cyl}(f)$. Moreover, using Plancherel identity, we have that
    \begin{align*}
\|u\|_{\mathcal{H}^\beta_{cyl}(f)}^2 &= \int_{\mathbb{R}^n} \left| (1+f(-|\xi|^2))\mc{F}(u)(\xi) \right|^2 d\xi 
= \int_{\mathbb{R}^n} \left| \mc{F}(g)(\xi) \right|^2 d\xi 
= \|g\|_{L^2_{cyl}(\mathbb{R}^n)}^2.
\end{align*}

\end{proof}

\begin{theorem}\label{teorema para Rn con simetria cilindrica}
   Suppose that $0 < \gamma < 1$ is a constant, $f \in \mathcal{G}^\beta$ for some $\beta > \frac{n\gamma}{2}$, and that $U:\mathbb{R}^{n} \times  \mathbb{R} 
\to \mathbb{R}$ is cylindrically symmetric in $x$, is $C^{1}$ in the second variable and $U(\cdot, 0) \not\equiv 0$.  
Assume moreover that there exist functions $h \in L^2(\mathbb{R}^n)$, $g \in L^{2/\gamma}(\mathbb{R}^n)$ and some constant $C>0$ such that the following two inequalities hold,
\begin{equation}\label{desigualdades caso cilindrico}
 |U(x, y) + y| \leq C(|h(x)| + |y|^{\frac{1}{1-\gamma} } ),
\end{equation}
\begin{equation}\label{desigualdades caso cilindrico 2}
    \left|\frac{\partial}{\partial y} (U(x, y) + y)\right| \leq C(|g(x)| + |y|^{\frac{\gamma}{1 - \gamma} }).
\end{equation}
Then there exist $0 < \tau < 1$ and $0 < \rho_{\tau } < 1$ such that, whenever $\|h\|_{L^2(\mathbb{R}^n)} < \rho_{\tau}$, there exists a solution $u \in \mathcal{H}^\beta_{cyl}(f)$ to the equation $f(\Delta)u - U(\cdot, u) = 0$ with $\|u\|_{L^{2/(1-\gamma) }_{cyl}(\mathbb{R}^n)} \leq \tau $.
\end{theorem}

\begin{proof}
Consider equation \eqref{cylindrical eq} and the operator $L=f(\Delta) + I$ defined in (\ref{operador no local L}). Replacing the function $Z_{u}(x)=U(x,u(x)) + u(x)$ with $u \in L^{2/(1-\gamma)}_{cyl}(\mathbb{R}^{n})$, the nonlinear equation \eqref{cylindrical eq} is equivalent to the  problem
\begin{equation}\label{sistema lineal caso cilindrico}
 Lu=Z_{u}.   
\end{equation}
The function $Z_{u}$ is cylindrically symmetric because both $U$ and $u$ have this symmetry. Moreover $Z_{u} \in L^{2}(\mathbb{R}^{n})$, since
\begin{align*}
|Z_{u}|^{2} &= |U(x,u(x)) + u(x)|^{2} \leq C(|h(x)| + |u(x)|^{\frac{1}{(1-\gamma)}} )^{2} \leq \widetilde{C}(|h(x)|^{2} + |u(x)|^{\frac{2}{(1-\gamma)}} )
\end{align*}
with $\widetilde{C}> 0$. Since $u \in L^{2/(1-\gamma)}_{cyl}(\mathbb{R}^{n})$ we have $|u|^{1/(1-\gamma)} \in L^{2}(\mathbb{R}^{n})$, also we have that $h(x) \in L^{2}(\mathbb{R}^{n})$, hence $Z_{u} \in L^{2}(\mathbb{R}^{n})$. We shall apply Schauder's fixed point theorem for compact operators, and our main tool will be Lion's theorem on Sobolev embeddings for cylindrically symmetric functions \cite[\textsc{Theorem III.3}]{Lions}. Given $\tau > 0$, define
\begin{equation}
    \mathcal{C}_{\tau}= \{ u \in L^{2/(1-\gamma) }_{cyl}(\mathbb{R}^{n}) \mid \norm{u}_{L^{2/(1-\gamma)}(\mathbb{R}^{n})} \leq \tau \}
\end{equation}
This is a nonempty subset of the Banach space \textcolor{black}{$L^{2/(1-\gamma)}(\mathbb{R}^{n})$} which is closed, bounded and convex. Define the map
\begin{equation}
    \mathcal{Q} : \mathcal{C}_{\tau} \to L^{2/(1-\gamma) }_{cyl}(\mathbb{R}^{n}), \quad \mathcal{Q}(u)=\widetilde{u}
\end{equation}
where $\widetilde{u}$ is the unique solution of the linear non-homogeneous problem $L \widetilde{u} = Z_{u}$.
Since $Z_{u} \in L^{2}_{cyl}(\mathbb{R}^{n})$ and by Theorem \ref{teorema de cerradura}, $\mathcal{Q}$ is well defined. We now show there exists $\tau >0$ such that \textcolor{black}{$\norm{\mathcal{Q}(u)}_{L^{2/(1-\gamma) }(\mathbb{R}^{n})} \leq \tau $}. Since $2 < \frac{2}{1- \gamma} < 2^{*}=\frac{2n}{n-2s_{\gamma}}$, in which $s_{\gamma} \in \left( \frac{n \gamma}{2}, \beta  \right)$, the Sobolev embedding theorem(\cite[p.5]{Taylor III})
and using Lemma (\ref{Lemma de contencion epscios Hb}) yields the following embedding chain
\begin{equation}\label{embedding 1}
   \mathcal{H}^{\beta}(f) \hookrightarrow H^{s_{\gamma}}(\mathbb{R}^{n}) \hookrightarrow
    H^{\frac{n \gamma}{2}}(\mathbb{R}^{n})  \hookrightarrow L^{2/(1-\gamma) }(\mathbb{R}^{n})
\end{equation}
Taking $u \in \mathcal{C}_{\tau}$,
\begin{equation}
    \norm{ \mathcal{Q}(u)}_{L^{2/(1-\gamma)}(\mathbb{R}^{n})} \leq c \norm{\mathcal{Q}(u)}_{H^{s_{\gamma}}(\mathbb{R}^{n})} \leq c_{\gamma}\norm{\mathcal{Q}(u)}_{\mathcal{H}^{\beta}_{cyl}(f)}
\end{equation}
where $c_{\gamma}$ is a constant coming from the embedding (\ref{embedding 1}). Using that $\mathcal{Q}(u)=\widetilde{u}$ is the unique solution of the linear problem, by Theorem \ref{teorema caso lineal} and inequality (\ref{desigualdades caso cilindrico}) we obtain
\begin{align*}
    \norm{\mathcal{Q}(u)}_{L^{2/(1-\gamma) }(\mathbb{R}^{n})} &\leq   c_{\gamma} \norm{\mathcal{Q}(u)}_{\mathcal{H}^{\beta}_{cyl}(f)} = c_{\gamma} \norm{Z_{u}}_{L^{2 }(\mathbb{R}^{n})}=c_{\gamma}\norm{U(\cdot , u(\cdot)) + u(\cdot)}_{L^{2}(\mathbb{R}^{n})} \\
    & \leq c_{\gamma} K \left( \norm{h}_{L^{2}(\mathbb{R}^{n})} + \norm{u}_{L^{2/(1-\gamma) }(\mathbb{R}^{n})}^{\frac{1}{(1-\gamma)} }   \right) \leq c_{\gamma} K \left( \norm{h}_{L^{2}(\mathbb{R}^{n})} + \tau^{1/(1-\gamma) }  \right) 
\end{align*}
for some constant $K>0$. If $\norm{h}_{L^{2}(\mathbb{R}^{n})}$ is sufficiently small, we may choose $0 < \tau < 1$ and $0 < \rho_{\tau} < 1$ such that $\norm{h}_{L^{2}(\mathbb{R}^{n})} < \rho_{\tau}$ and
$$c_{\gamma} K \left(\norm{h}_{L^{2}(\mathbb{R}^{n})} + \tau^{1/(1-\gamma) } \right) < \tau.$$
Hence $\norm{\mathcal{Q}(u)}_{L^{2/(1-\gamma) }(\mathbb{R}^{n})} \leq \tau$, so $\mathcal{Q}: \mathcal{C}_{\tau} \to \mathcal{C}_{\tau}$ is well defined. \\
We now show $\mathcal{Q}$ is compact. Let $(u_{n}) \subset \mathcal{C}_{\tau} \subset L^{2/(1- \gamma) }_{cyl}(\mathbb{R}^{n})$ be a bounded sequence, since the sequence $\mathcal{Q}(u_{n})$ is bounded in $L^{2/(1-\gamma)}(\mathbb{R}^{n})$, we have $ c_{\gamma}\| \mc{Q}(u_{n})\|_{\mc{H}^{\beta}_{cyl}(f)} \leq c_{\gamma}K(\|h\|_{L^{2}} + \tau^{1/(1-\gamma)})$
and hence bounded in $H^{s_{\gamma}}_{cyl}(\mathbb{R}^{n})$ by (\ref{embedding 1}). By \cite[\textsc{Theorem III.3}]{Lions}, the embedding \begin{equation}
  H^{s_{\gamma}}_{cyl}(\mathbb{R}^{n}) \hookrightarrow L^{2/(1-\gamma) }_{cyl}(\mathbb{R}^{n})  
\end{equation}
is compact, so there exists a subsequence $\mathcal{Q}(u_{n_{k}})$ 
of $\mathcal{Q}(u_{n})$ that converges in $L_{cyl}^{2/(1-\gamma)}(\mathbb{R}^{n})$. 
Therefore $\mathcal{Q}$ is compact. 
Finally, we prove continuity of $\mathcal{Q}$. Let $(u_{n}) \subset \mathcal{C}_{\tau}$ with $u_{n} \to u$ in $L^{2/(1-\gamma)}_{cyl}(\mathbb{R}^{n})$.
Our main tool is the embedding chain (\ref{embedding 1}). Then
\begin{equation*}
    \norm{\mathcal{Q}(u_{n}) - \mathcal{Q}(u)}^{2}_{L^{\frac{2}{(1-\gamma)}}(\mathbb{R}^{n})} \leq c \norm{\widetilde{u_{n}} - \widetilde{u}}^{2}_{\mathcal{H}^{\beta}_{cyl}(f)}=c\norm{Z_{u_{n}} - Z_{u}}^{2}_{L^{2}(\mathbb{R}^{n})}, \quad c>0~const. 
\end{equation*}
the last equality is because the linearity of the system $Lu=Z_{u}$. Because if 
$\widetilde{u}_{n}, \widetilde{u}$ are the uniquest solutions of $Lu=Z_{u_{n}}$ and $Lu=Z_{u}$ respectively, we have $L(\widetilde{u}_{n} - \widetilde{u})=Z_{(u_{n}-u)}$, and $Z_{(u_{n}-u)}=L\widetilde{u}_{n}-L\widetilde{u}=Z_{u_{n}}-Z_{u}$. Denote $W(x,u(x))=U(x,u(x)) + u(x)$, then we have 
\begin{align*}
    \norm{\mathcal{Q}(u_{n}) - \mathcal{Q}(u)}^{2}_{L^{\frac{2}{(1-\gamma)}}(\mathbb{R}^{n})} &\leq c\norm{Z_{u_{n}} - Z_{u}}^{2}_{L^{2}(\mathbb{R}^{n})}= c\int_{\mathbb{R}^{n}}|W(x,u_{n}(x)) - W(x,u(x)) |^{2}dx \\
    &= c \int_{\mathbb{R}^{n}}\left| \int_{0}^{1}\frac{d}{dt}W(x,tu_{n}(x) 
   +  (1-t)u(x))dt  \right|^{2}dx \\
   &= c \int_{\mathbb{R}^{n}} |u_{n}(x) - u(x)|^{2}\left| \int_{0}^{1}\frac{\partial}{\partial y}W(x, tu_{n}(x) 
   +  (1-t)u(x))    dt \right|^{2}dx
\end{align*}
We define the function $j_{x}: [0,1] \to \mathbb{R}, j_{x}(t)=\frac{\partial}{\partial y}W(x, tu_{n}(x) 
   +  (1-t)u(x))$ with $x \in 
\mathbb{R}^{n}$, since $U(\cdot, u)$ is $C^{1}$ in the second variable, $j_{x}(t)$ is a well defined Carathéodory function(i.e., measurable with respect to $x$ and continuous with respect to $t$). Note that
\begin{equation*}
    \Big|\int_{0}^{1}j_{x}(t)dt \Big| \leq \int_{0}^{1} |j_{x}(t)|dt \leq \int_{0}^{1}\left( \sup_{t \in [0,1]}|j_{x}(t)|  \right)dt = \sup_{t \in [0,1]}|j_{x}(t)|
\end{equation*}
Since $j_{x}(t)$ is continuous in a compact set, the maximum is attained for each $x \in \mathbb{R}^{n}$, i.e., for each fixed $x_{0} \in \mathbb{R}^{n}$,
there exists $t_{0} \in [0,1]$ such that $j_{x_{0}}$ attained the supremum in $t_{0}$.
In fact there exists a Borel measurable function $t(\cdot): M \to [0,1]$ by \cite[t.18.13]{AliprantisBorder} 
such that for each $x^{\prime} \in \mathbb{R}^{n}$, $t(x^{\prime})$ 
satisfied $\sup_{t \in [0,1]}|j_{x^{\prime}}(t)|=|j_{x^{\prime}}(t(x^{\prime}))|$. Using (\ref{desigualdades caso cilindrico 2}) for a constant $C > 0$ we get
{\fontsize{11}{14}\selectfont
\begin{align*}
  \norm{\mathcal{Q}(u_{n}) - \mathcal{Q}(u)}^{2}_{L^{\frac{2}{(1-\gamma)}}(\mathbb{R}^{n})}   &\leq 
C\int_{\mathbb{R}^{n}}|u_{n}(x) - u(x)|^{2} \left\{  \sup_{t \in [0,1]}\left|  |g(x)| + |tu_{n}(x) + (1-t)u(x)|^{\frac{\gamma}{1-\gamma} } \right| \right\}^{2}dx\\
  &\leq C\int_{\mathbb{R}^{n}}|u_{n}(x) - u(x)|^{2} \left\{   |g(x)| + |t(x)u_{n}(x) + (1-t(x))u(x)|^{\frac{\gamma}{1-\gamma} } \right\}^{2}dx
\end{align*}}
For the second inequality, we use the fact $0 \leq t(x) \leq 1$ to obtain the following computations
\begin{align*}
|g(x)| + |t(x)u_{n}(x) + (1-t(x))u(x)|^{\frac{\gamma}{1-\gamma}} &\leq |g(x)| + \left( |t(x)u_{n}(x)| + |(1-t(x))u(x)|  \right)^{\frac{\gamma}{(1-\gamma)}} \\
& \leq |g(x)| + c_{1}\left( 
|t(x)u_{n}(x)|^{\frac{\gamma}{(1-\gamma)}} + |(1-t(x))u(x)|^{\frac{\gamma}{(1-\gamma)}}
\right) \\
&\leq |g(x)| + c_{2}\left( 
|u_{n}(x)|^{\frac{\gamma}{(1-\gamma)}} + |u(x)|^{\frac{\gamma}{(1-\gamma)}}
\right),
\end{align*}
with $c_{1}, c_{2}>0$ constants. Then for a constant $\widetilde{c}>0$, using H\"older inequality, we have 
\begin{align*}
  \norm{\mathcal{Q}(u_{n}) - \mathcal{Q}(u)}^{2}_{L^{\frac{2}{(1-\gamma)}}(\mathbb{R}^{n})}&\leq \widetilde{c}\int_{\mathbb{R}^{n}}|u_{n}(x) - u(x)|^{2}\left( |g(x)|^{2} + |u_{n}(x)|^{\frac{2\gamma}{1-\gamma}} + |u(x)|^{\frac{2\gamma}{1-\gamma}}\right) dx  \\
   &= \widetilde{c}\Big(\int_{\mathbb{R}^{n}}|u_{n}(x) - u(x)|^{2}|g(x)|^{2} dx + \int_{\mathbb{R}^{n}}|u_{n}(x) - u(x)|^{2}|u_{n}(x)|^{\frac{2\gamma}{1-\gamma} } dx \\
   &+ \int_{\mathbb{R}^{n}}|u_{n}(x) - u(x)|^{2}|u(x)|^{\frac{2\gamma}{1-\gamma}} dx \Big) \\ 
   & \leq \widetilde{c} \norm{u_{n} - u}^{2}_{L^{2/(1-\gamma)}(\mathbb{R}^{n})} \Big(\norm{g}^{2}_{L^{2/\gamma}(\mathbb{R}^{n})}  + \tau^{\frac{2\gamma}{1-\gamma}}    \Big) .
\end{align*}
Hence $\mathcal{Q}(u_{n}) \to \mathcal{Q}(u)$ in $L^{2/(1-\gamma)}(\mathbb{R}^{n})$.
Therefore $\mathcal{Q}: \mathcal{C}_{\tau} \to \mathcal{C}_{\tau}$ satisfies the hypotheses of Schauder's fixed point theorem. 
\end{proof}

\newpage

\section{Compact Riemannian Manifolds with Symmetries}\label{seccion de compact Riemannian manifolds}

Let $(M,g)$ be a compact, connected and oriented Riemannian $n$-manifold without boundary. This section extends the theory of nonlocal equations with symmetry to the setting of Riemannian manifolds,
incorporating group actions of subgroups $G \leq \text{Isom}_{g}(M)$ over $(M,g)$. Here $\text{Isom}_g(M)$ is the group of all Riemannian isometries of $M$. 
We consider the
equation
\begin{equation}
f(\Delta_g)\phi - U(\cdot, \phi) = 0
\end{equation}
Where the operator $f(\Delta_g)$ is defined through the spectral 
theorem and $U: M \times \mathbb{R} \to \mathbb{R}$ represents a nonlinearity which is $C^{1}$ in the second variable and satisfying 
specific growth conditions. Moreover we say that $U: M \times \mathbb{R} \to \mathbb{R}$ is $G$-invariant if $U(\sigma(x),t)=U(x,t)$ for every $\sigma\in G$, $x\in M$, $t\in\mathbb{R}$. 
Our main objective in this section is to obtain an existence result for compact Riemannian manifolds under the action of a group with certain 
properties, more precisely, under the action of closed subgroups $G\leq\mathrm{Isom}_{g}(M)$. 
As shown in \cite[p.249]{Hebey 2}, $\mathrm{Isom}_{g}(M)$ is a
Lie group when endowed with the compact open topology, and it acts differentiably on $M$. Since $(M,g)$ is a  compact manifold, $\mathrm{Isom}_{g}(M)$ is a compact Lie group, hence $G$ is itself a Lie group. 
The Sobolev spaces under actions of a group were defined formally in \cite{Hebey 2} as 
  \begin{equation}
        H^{s,p}_{G}(M)=\{u \in H^{s,p}(M) \mid u \circ \sigma = u, \forall \sigma \in G\}, \quad s \in \mathbb{R}_{>0}.
    \end{equation}
The spaces $H^{s,p}(M)$ where already defined in \cite[p.689]{Banica} as 
\begin{equation}\label{espacios de sobolev en version laplaciano R1}
    H^{s,p}(M)=\{u \in L^{p}(M) \mid \exists h \in L^{p}(M), u=(I+\Delta_{g})^{-s/2}h \}, 
\end{equation}
endowed with the norm $\|u\|_{H^{s,p}(M)}=\| h \|_{L^p(M)}.$
We have that the space of $G$-invariant $C^{\infty}_{0}(M)$ functions denoted by $\mathcal{D}_{G}(M)$ satisfied $\mathcal{D}_{G}(M) \subset H^{1,p}_{G}(M)$,
and if $G$ is compact, one gets from the existence of the Haar measure that $\mathcal{D}_{G}(M)$ is dense in $H^{1,p}_{G}(M)$.
The following lemmas will serve as tools to prove the main theorem of this section. The 
first of these study the density of $\mathcal{D}_{G}(M)$ under the fractional Sobolev spaces under actions $H^{s,p}_{G}(M)$.
The second lemma is the generalization of the Theorem \cite[t.9.1]{Hebey 2} to 
fractional Sobolev spaces under the action of a closed subgroup $G \leq \text{Isom}_{g}(M)$.

\begin{lemma}\label{espacio denso en espacios de sobolev fraccionarios}
Let $(M,g)$ be a compact Riemannian manifold and let $G \leq \text{Isom}_{g}(M)$ be a closed subgroup. For any $s,p \in \mathbb{R}$ with $s>0$ and $p\geq 1$, the space $\mathcal{D}_{G}(M)$ is dense in $H^{s,p}_{G}(M)$.
\end{lemma}

\begin{proof}
Note that $\mathcal{D}(M)$ is dense in $H^{s,p}(M)$ because \cite[t.2.6]{Aubin 2} proves density of $\mathcal{D}(M)$ for spaces $H^{k,p}(M)$ 
with $k \in \mathbb{Z}$ and any real $p \geq 1$, which implies that the same result holds for fractional Sobolev spaces $H^{s,p}(M)$ defined by complex interpolation method \cite[p.25]{Taylor III}. 
Fix $s \in \mathbb{R}_{>0}$ and a $p \geq 1$ real.\\
Let $u \in H^{s,p}_{G}(M)$ and $(u_{m})\subset \mathcal{D}(M)$ 
be a sequence such that $u_{m} \to u$ in $H^{s,p}(M)$. Since $M$ is compact, $\text{Isom}_{g}
(M)$ is a compact Lie group, then $G$ is a compact Lie 
subgroup, which implies the existence of a Haar measure $\sigma_{G}$ on $G$\cite[t.9.2.2]{Cohn}.
Define the averaging operator $\mathcal{A}_{G}$, which is a linear and continuous operator
\begin{equation}
    (\mathcal{A}_{G}u_{m})(x)=\frac{1}{\sigma_{G}(G)}\int_{G}u_{m}(\sigma(x))d\sigma_{G}, \quad x\in M.
\end{equation}
Consider $(G,\mathcal{B}_{G}, \sigma_{G})$, the measure space induced by $G$, with $\mathcal{B}_{G}$ the Borel $\sigma$-algebra.
We want to prove that $\mathcal{A}_{G}u_{m} \in \mathcal{D}_{G}(M)$ for each $m  \in \mathbb{N}$. To prove differentiability define $\xi: M \times G  \to \mathbb{R}, (x,\sigma) \mapsto \xi(x,\sigma)=u_{m}(\sigma(x))$. 
Note that for any $\sigma \in G$, the function $x \mapsto \xi(x,\sigma)$ is $C^{k}(M)$ for 
$k \geq 1$, because $x \mapsto u_{m}(\sigma(x))$ is $C^{k}(M)$ since $u_{m} \in C^{\infty}(M)$ and the action $G 
\curvearrowright M$ is smooth \cite[p.249]{Hebey 2}, i.e. $(x,\sigma) \mapsto \sigma(x)$ is smooth. If we 
consider a multi-index $\alpha \in \mathbb{N}_{0}^{n}$ such that $|\alpha| \leq k$ we 
have that the application $\psi: G \to \mathbb{R}$ defined as $\sigma \mapsto \partial^{\alpha}_{x}u_{m}(\sigma(x))$ is smooth, which implies that it is measurable in $(G,\mathcal{B}_{G}, \sigma_{G})$.
Using that $\xi$ is a smooth function in $M \times G$ and the compactness of $G$, the $\alpha$-derivatives in $x$ 
are uniformly bounded for $\sigma \in G$, then there exists $K>0$ such that 
$|\partial^{\alpha}_{x}\xi(x,\sigma)| \leq g_{\alpha}(\sigma)$ with $g_{\alpha}(\sigma)=K \in L^{1}(G,\sigma_{G})$. Therefore by Dominated Convergence Theorem we have
\begin{equation}
    \partial^{\alpha}_{x}(\mathcal{A}_{G}u_{m}(x))=\frac{1}{\sigma_{G}(G)}\int_{G}\partial^{\alpha}_{x}(u_{m}(\sigma(x)))d\sigma_{G}
\end{equation}
Since $u_{m}\in C^{\infty}(M)$ we have that $\mathcal{A}_{G}u_{m} \in C^{\infty}(M)$ for each $m \in \mathbb{N}$. Fix $m \in 
\mathbb{N}$, suppose that $\widetilde{K}_{m}=\text{supp} (\mathcal{A}_{G}u_{m})$ and $K_{m}=\text{supp}(u_{m})$. If $x \in \widetilde{K}_{m}$, then 
$\mathcal{A}_{G}u_{m}(x) \neq 0$, then there exists a $\sigma \in G$ such that 
$u_{m}(\sigma(x)) \neq 0$ and $ \sigma(x) \in K_{m}$, 
indeed $x \in \sigma^{-1}(K_{m})$. Moreover we have
\begin{equation}
    \widetilde{K}_{m} \subset \bigcup_{\sigma \in G}\sigma^{-1}(K_{m}).
\end{equation}
Define $\vartheta: G \times K_{m} \to M$ by $(\sigma, y) \mapsto 
\sigma^{-1}(y)$, this function is continuous; indeed since $G \times K_{m}$ is compact, we have that $\vartheta(G\times K_{m})=\cup_{\sigma \in G }\sigma^{-1}(K_{m})$ is compact\footnote{This argument can be useful in the case $(M,g)$ is non-compact.}.
Since $\widetilde{K}_{m}$ is a closed subset of a compact set, $\widetilde{K}_{m}$ is compact on $M$. Therefore 
$\mathcal{A}_Gu_m \in \mathcal{D}(M)$ for each $m \in \mathbb{N}$. Moreover for every $\tau \in G$, we have the following computations 
\begin{align*}
    (\mathcal{A}_{G}u_{m})(\tau(x))&= \frac{1}{\sigma_{G}(G)}\int_{G}u_{m}(\sigma(\tau(x)))d\sigma_{G}=  \frac{1}{\sigma_{G}(G)}\int_{G}u_{m}(\sigma^{\prime}(x))d\sigma_{G}=(
    \mathcal{A}_{G}u_{m}
    )(x).
\end{align*}
Hence $\mathcal{A}_{G}u_{m} \in \mathcal{D}_{G}(M)$ for any $m \in \mathbb{N}$, then $(\mathcal{A}_{G}u_{m})_{m \in \mathbb{N}} \subset 
\mathcal{D}_{G}(M)$. We show that $\mathcal{A}_{G}u_{m} \to u$ in $H^{s,p}(M)$. Note that if $u \in H^{s,p}_{G}(M)$ then $(\mathcal{A}_{G}u)(x)=u(x)$. Hence the following 
\begin{align*}
\norm{\mathcal{A}_{G}u_{m} - u}_{H^{s,p}(M)} &\leq \norm{\mathcal{A}_{G}u_{m} - u_{m}}_{H^{s,p}(M)} + \norm{u_{m}-u}_{H^{s,p}(M)}\\
&\leq \norm{\mathcal{A}_{G}u_{m} - \mathcal{A}_{G}u}_{H^{s,p}(M)} + \norm{\mathcal{A}_{G}u - u_{m}}_{H^{s,p}(M)} + \norm{u_{m} - u}_{H^{s,p}(M)}\\
&=\norm{\mathcal{A}_{G}(u_{m}-u)}_{H^{s,p}(M)} + 2\norm{u_{m}-u}_{H^{s,p}(M)} \leq 3\norm{u_{m}-u}_{H^{s,p}(M)}.
\end{align*}
Since $u_{m} \to u$ in $H^{s,p}(M)$, it follows that $\mathcal{A}_{G}u_{m} \to u$ in $H^{s,p}(M)$. 
\end{proof}

\begin{lemma}\cite[l.9.1]{Hebey 2}.\label{lemma de las cartas} 
Let $(M,g)$ be a smooth Riemannian $n$-manifold (not necessarily compact or complete), and let $G$ be a compact subgroup of $\mathrm{Isom}_g(M)$. Let $x \in M$ and 
$k=\dim \mathcal{O}_{G}(x)$. Assume $k \geq 1$. There exists a coordinate chart $(\Omega, \varphi)$ of $M$ at $x$ such that: \\
$i)$ $\varphi(\Omega)=U \times V$, where $U$ is some open subset of $\mathbb{R}^{k}$ and $V$ is some open subset of $\mathbb{R}^{n-k}$, \\
$ii)$ $\forall y \in \Omega, U \times \pi_{2}(\varphi(y)) \subset \varphi(\mathcal{O}_{G}(y) \cap \Omega)$ where $\pi_{2}: \mathbb{R}^{k} \times \mathbb{R}^{n-k} \to \mathbb{R}^{n-k}$ is the second projection.
\end{lemma}

\justify

Let $\textbf{X}=(X_{0},X_{1})$ and $\textbf{Y}=(Y_{0},Y_{1})$ be Banach couples, i.e., pairs of Banach spaces continuously embedded into a common Hausdorff locally convex topological vector space. Let $T: \textbf{X} \to \textbf{Y}$ be a bounded linear operator in the sense of Cwikel, meaning that $T: X_{0} + X_{1} \to Y_{0} + Y_{1}$ is a linear map such that its restrictions $T|_{X_j}: X_{j} \to Y_{j}$ are bounded operators for $j=0,1$. As is standard, we denote the complex interpolation space by $$X_{\theta}=[X_{0},X_{1}]_{\theta}, \quad \theta \in (0,1).$$

\begin{theorem}\cite[t.10]{Cwikel}.\label{Cwikel theorem for interpolation compacity1}
    Let $\textbf{X}$ be a Banach couple such that $X_{0}$ is reflexive and is given by $X_{0}=[W,X_{1}]_{\alpha}$ for some $0 < \alpha < 1$ and some Banach space $W$ which forms a Banach couple with $X_{1}$. Suppose $T: \textbf{X} \to \textbf{Y}$ is such that $T: X_{0} \to Y_{0}$ is compact. Then 
    \begin{equation}
       T: X_{\theta} \to Y_{\theta} 
    \end{equation}
     is compact for $0<\theta < 1$.
\end{theorem}

To establish our main existence result on compact Riemannian manifolds with symmetries (Theorem~\ref{Teorema Final}) via Schauder's fixed-point theorem, we require a compactness mechanism adapted to the nonlocal operator $f(\Delta_g)$ and the $G$-invariant functional setting.

Since the classical Rellich--Kondrachov theorem does not exploit the dimensional gain produced by infinite $G$-orbits, we first prove the following fractional $G$-invariant compact embedding by combining the integer-order Hebey--Vaugon theorem~\cite[t.9.1]{Hebey 2} with complex interpolation theory.
This theorem serves as the fundamental analytical tool for the proof of Theorem~\ref{Teorema Final}:

\newpage  

\begin{theorem}\label{Lema principal caso con accion}

 Let $(M,g)$ be a compact, connected, oriented Riemannian $n$-manifold without boundary. Let $q>1$, $s>0$ be real numbers, and let $G$ be a closed subgroup of $\mathrm{Isom}_g(M)$. Assume that for any $x \in M$, $|\mathcal{O}_{G}(x)|=+\infty$ and set $k=\min_{x\in M}\dim \mathcal{O}_{G}(x).$ Then, $k\ge 1$ and for every $p \in [1,\infty)$ satisfying
$
s>(n-k)\left(\frac{1}{q}-\frac{1}{p}\right),$
the embedding
\begin{equation}
   H^{s,q}_{G}(M)\hookrightarrow L^{p}(M) 
\end{equation}
is compact. Equivalently, if $qs \geq n-k$, the embedding is compact for every real $p \in [1, 
\infty)$. More precisely, if $qs > n-k$, $H^{s,q}_{G}(M)$ embeds compactly into $C^{0}(M)$(and hence into $L^{\infty}(M)$).
If $qs<n-k$, the embedding is compact for $1 \leq p < \frac{(n-k)q}{n-k-sq}$, and continuous for $1 \leq p \leq \frac{(n-k)q}{n-k-sq}$. 

\end{theorem}

\begin{proof}
We proceed as follows. We first generalize the result \cite[t.9.1]{Hebey 2}, which proves that the embedding
$H^{1,q}_{G}(M)\hookrightarrow L^{p}(M)$
is compact under suitable conditions on $p,q\geq 1$ real numbers, to Sobolev spaces of integer order $\ell\ge 1$. Using this generalization and the complex interpolation characterization of fractional Sobolev spaces \cite[p.25]{Taylor III},
$$
H^{s,q}(M)=[H^{\kappa,q}(M),H^{\kappa+1,q}(M)]_{\theta},
$$
where $\kappa\ge 0$ is an integer, $q>1$, $\theta\in(0,1)$, and $s=\kappa+\theta$, $s>0$, we extend the integer-order embeddings to the fractional-order case.
For the endpoint case $p=\infty$, we use a Hölder space argument: when $qs>n-k$, the fractional Sobolev space $H^{s,q}_{G}(M)$ embeds compactly into $C^{0}(M)$, and hence into $L^{\infty}(M)$.
Note that $k \geq 1$.
In fact, suppose that there exists $x_{0} \in M$ such that $\dim \mathcal{O}_{G}(x_{0})=0$. Since $G$ is a closed subgroup of the compact Lie group $\mathrm{Isom}_g(M)$, $G$ is compact, so $\mathcal{O}_{G}(x_{0})$ is a 
compact $0$-dimensional manifold, hence a finite discrete manifold, which 
contradicts the fact that $|\mathcal{O}_{G}(x)|=\infty$ for every $x \in M$.
If $k=n$, every orbit has dimension $n$ and is a smooth open and compact submanifold of $M$; hence it is closed. Therefore, by the connectedness of the manifold, the action $G\curvearrowright M$ is transitive, which implies that for each $p \in M$, $\mc{O}_{G}(p)=M$.
Henceforth suppose $1\leq k<n$.\\
Since $M$ is compact, using Lemma \eqref{lemma de las cartas} there exists a finite covering $\{(\Omega_{m}, \varphi_{m})\}_{m=1}^{N}$, such that each $(\Omega_{m}, \varphi_{m})$ satisfies the following conditions \\
$a)$ $\varphi_{m}(\Omega_{m})=U_{m} \times V_{m}$, where $U_{m} \subset \mathbb{R}^{k_{m}}$ is open and bounded, and $V_{m} \subset \mathbb{R}^{n-k_{m}}$ is open and bounded with smooth boundary, where $k_{m} \in \mathbb{Z}$, $k \leq k_{m} \leq n-1$,\\
$b)$ $\forall z 
\in \Omega_{m}, U_{m} \times \pi_{2}(\varphi_{m}(z)) \subset \varphi_{m}(\mathcal{O}_{G}(z) \cap \Omega_{m})$, \\
$c)$ $\exists \alpha_{m} > 0$ with $\alpha_{m}^{-1}\delta_{ij} \leq g_{ij}^{m} \leq \alpha
_{m}\delta_{ij}$ as bilinear forms, where $g_{ij}^{m}$ are the metric coefficients in each chart $(\Omega_{m}, \varphi_{m})$. Without loss of generality, one may assume that $\varphi_{m}$ is defined on some open $\widetilde{\Omega}_{m}$ containing $\overline{\Omega}_{m}$ such that $\varphi_{m}(\widetilde{\Omega}_{m})=\widetilde{U}_{m}\times \widetilde{V}_{m}$ with  $\overline{U}_{m} \subset \widetilde{U}_{m}$ and $\overline{V}_{m} \subset \widetilde{V}_{m}$.

\vspace{0.12cm}
$i)$\label{punto 1, sobolev entero con simetria 1} Let $u \in \mathcal{D}_{G}(M)$. If we take $x,x^{\prime} \in U_{m}$, using $b)$ we have that for every $z \in \Omega_{m}$, $\varphi_{m}^{-1}(U_{m} \times \pi_{2}(\varphi_{m}(z))) \subset \{\sigma(z) \in \Omega_{m} \mid \sigma \in G \}$.
Indeed, for any $y \in V_{m}$ there exists $h \in \Omega_{m}$ such that $y=\pi_{2}(\varphi_{m}(h))$, and for every $x \in U_{m}$ one has $\varphi_{m}^{-1}(x,y)=\varphi_{m}^{-1}(x,\pi_{2}(\varphi_{m}(h))) \in \{ \sigma(h) \in \Omega_{m} \mid \sigma \in G\}$. 
Then there exist $\sigma_{1}, \sigma_{2} \in G$ such that $\varphi_{m}^{-1}(x,y)=\sigma_{1}(h)$ and $\varphi_{m}^{-1}(x^{\prime},y)=\sigma_{2}(h)$, hence by the $G$-invariance of $u$,
$$u \circ \varphi_{m}^{-1}(x,y) = u \circ \varphi_{m}^{-1}(x^{\prime},y), \quad m \in \{1,...,N\}. $$
Thus, there exists $\widetilde{u}_{m} \in C^{\infty}(\mathbb{R}^{n-k_{m}})$ defined by 
\begin{equation}\label{funcion truco invariante caso riemaniano simetrico}
  \widetilde{u}_{m}(y)=u \circ \varphi_{m}^{-1}(x,y).  
\end{equation}
Fix any $m \in \{1,...,N\}$ and any $p \geq 1$ real.
Take local coordinates on $M$ and on $\mathbb{R}^{n}$ given by $(y^{i})$ and $(x^{j})$, respectively.
Consider the diffeomorphism $\varphi_{m}^{-1}: \varphi_{m}(\Omega_{m}) \to \Omega_{m}$. By \cite[p.14.20]{Lee} and condition $c)$, which guarantees the existence of $\alpha_{m}>0$ such that $\alpha_{m}^{-1}\delta_{ij} \leq g_{ij}^{m} \leq \alpha
_{m}\delta_{ij}$, it follows that $\sqrt{\text{det}\,g_{ij}^{m}}$ is bounded from above and below; denoting $\omega=|u|^{p}dV_{g}$ we obtain  
\begin{align*}
\int_{\Omega_{m}}|u|^{p}dV_{g} &= \int_{\varphi_{m}(\Omega_{m})}(\varphi_{m}^{-1})^{*}\omega=\int_{\varphi_{m}(\Omega_{m})}(|u|^{p} \circ  \varphi_{m}^{-1}(x,y))\sqrt{\det g_{ij}^{m}}dxdy \\
&\leq A_{m}\int_{U_{m} \times V_{m}}|u \circ \varphi_{m}^{-1}(x,y)|^{p}dxdy= A_{m}\int_{U_{m} \times V_{m}}|\widetilde{u}_{m}(y)|^{p}dxdy \\ 
&=A_{m}\int_{V_{m}}\left( \int_{U_{m}}|\widetilde{u}_{m}(y)|^{p}dx\right)dy = A_{m}vol_{\overline{g}}(U_{m})\int_{V_{m}}|\widetilde{u}_{m}(y)|^{p}dy      \leq \widetilde{A}_{m}\int_{V_{m}}|\widetilde{u}_{m}(y)|^{p}dy
\end{align*}
for some constants $A_{m},\widetilde{A}_{m}>0$. For higher derivatives we compare the full Sobolev norms, where $\nabla^{j}u$ denotes the $(0,j)$-tensor total covariant derivative of $u$. By condition $c)$ and the boundedness of the Christoffel symbols and their derivatives on $\overline{\Omega}_{m} \subset \widetilde{\Omega}_{m}$, the coefficients of these formulas and their inverses are bounded on $\overline{\Omega}_{m}$. Since $u\circ\varphi_{m}^{-1}$ is independent of $x$, its coordinate Sobolev norm on $U_{m}\times V_{m}$ is equivalent to that of $\widetilde{u}_{m}$ on $V_{m}$.
Following the same procedure as above we obtain for every integer $\ell\geq 0$ and any real $p \geq 1$
the existence of constants $\{A_{m,0},...,A_{m,\ell}, B_{m,0},...,B_{m,\ell}\}$, with $A_{m,j}, B_{m,j}>0$, $\forall j \in \{0,..,\ell\}$, such that 
\begin{align}\label{inequacion caso entero lema embeding accion}
\sum_{j=0}^{\ell}B_{m,j}\left(\int_{V_{m}}|\nabla^{j}\widetilde{u}_{m}(y)|^{p}dy\right)^{\frac{1}{p}} \leq  \sum_{j=0}^{\ell}\left( \int_{\Omega_{m}}|\nabla^{j}u|^{p}dV_{g} \right)^{\frac{1}{p}} \leq \sum_{j=0}^{\ell}A_{m,j}\left(\int_{V_{m}}|\nabla^{j}\widetilde{u}_{m}(y)|^{p}dy\right)^{\frac{1}{p}}.    
\end{align}
\justify
Using the Rellich--Kondrachov theorem for bounded domains in Euclidean spaces, the embedding $H^{\ell,q}(V_{m})  \hookrightarrow H^{0,p}(V_{m})$ is compact if and only if $\ell - \frac{n-k_{m}}{q} > -\frac{n-k_{m}}{p}$\,(and continuous if $\ell - \frac{n-k_m}{q} \geq -\frac{n-k_m}{p}$). We obtain  
\begin{align*}
    \left(\int_{\Omega_{m}}|u|^{p}dV_{g} \right)^{\frac{1}{p}} &\leq  K_{m}\left(\int_{V_{m}}|\widetilde{u}_{m}|^{p}dy \right)^{\frac{1}{p}} \leq \widetilde{K}_{m}\left[ \left( \int_{V_{m}}|\widetilde{u}_{m}|^{q}dy  \right)^{\frac{1}{q}} + \ldots + \left(\int_{V_{m}}|\nabla^{\ell}\widetilde{u}_{m}|^{q}dy   \right)^{\frac{1}{q}} \right]\\
    &\leq \widetilde{K}_{m}\left[ \widetilde{B}_{m,0}\left( \int_{\Omega_{m}}|u|^{q}dV_{g}  \right)^{\frac{1}{q}}+  \ldots + \widetilde{B}_{m,\ell}\left( \int_{\Omega_{m}}|\nabla^{\ell}u|^{q}dV_{g}    \right)^{\frac{1}{q}}   \right].
\end{align*}
Here $K_{m}=(\widetilde{A}_{m})^{\frac{1}{p}}$, $\widetilde{B}_{m,i}=(B_{m,i})^{-1}$ for every $i \in \{0,..,\ell\}$ and $\widetilde{K}_{m} > 0$ is a constant. 
Thus, for each $m\in\{1,\ldots,N\}$ and $q>1$, the following local estimates hold. \\
If $q \ell \geq n-k_{m}$, then for any real number $p\geq 1$ there exists a constant $C_{m}>0$ such that, for any $u \in \mathcal{D}_{G}(M)$, one has 
\begin{equation}\label{embedding para prueba de lema final 1}
    \norm{u}_{L^{p}(\Omega_{m})} \leq C_{m}\norm{u}_{H^{\ell,q}(\Omega_{m})}.
\end{equation}
If $q\ell < n-k_{m}$, then for any real number $p \geq 1$ satisfying $1 \leq p \leq \frac{(n-k_{m})q}{n-k_{m}-\ell q}$, there exists a constant $E_{m} > 0$ such that for any $u \in \mathcal{D}_{G}(M)$ one has 
\begin{equation}\label{embedding para prueba de lema final 2}
       \norm{u}_{L^{p}(\Omega_{m})} \leq E_{m}\norm{u}_{H^{\ell,q}(\Omega_{m})}. 
\end{equation}
Note that if $\ell q < n-k_{m} < n-k$, then $\frac{(n-k_{m})q}{n-k_{m}-\ell q} > \frac{(n-k)q}{n-k-\ell q}$ (whereas if $n-k_{m} \leq \ell q$, any $p \in [1,\infty)$ is admissible on $V_{m}$). Therefore, if $q\ell < n-k$ and $1 \leq p < \frac{(n-k)q}{n-k-\ell q}$(or if $q\ell \geq n-k$ and $p \in [1,\infty)$), the compactness 
conditions for the embedding $H^{\ell,q}(V_{m})  \hookrightarrow H^{0,p}(V_{m})$ are satisfied for every $m \in \{1,\dots,N\}$. 
We have that $\{(\Omega_{m},\varphi_{m})\}_{m=1}^{N}$ is a finite open covering of $\text{supp}\,(|u|^{p}dV_{g})$.
If we define a nonnegative smooth partition of unity $( \psi_{m} )_{m=1}^{N}$ subordinate
to this covering, we have the following computations with $\widetilde{C} = \max_{1\leq m \leq N}\max\{C_m, E_m\}$:
\begin{align*}
    \norm{u}_{L^{p}(M)} &= \left( \sum_{m=1}^{N} \int_{\Omega_{m}}\psi_{m}|u|^{p}\,dV_{g}  \right)^{\frac{1}{p}} \leq \sum_{m=1}^{N}\norm{u}_{L^{p}(\Omega_{m})}\leq  \widetilde{C}\sum_{m=1}^{N}\norm{u}_{H^{\ell,q}(\Omega_{m})} \leq \widetilde{C}N\norm{u}_{H^{
    \ell,q}(M)}.
\end{align*}
Therefore if $q\ell < n-k$, for any real number $p$ such that $1 \leq p < \frac{(n-k)q}{n-k-\ell q}$ and if $q\ell \geq n-k$ for any real $p \geq 1$, by the density of $\mc{D}_{G}(M)$ in $H^{\ell,q}_{G}(M)$ established by Lemma \ref{espacio denso en espacios de sobolev fraccionarios}, the embedding $H^{\ell,q}_G(M) \hookrightarrow L^p(M)$ is continuous (including at the critical exponent $p = \frac{(n-k)q}{n-k-\ell q}$ when $q\ell < n-k$). We now prove compactness under the strict inequality $1/p>1/q-\ell/(n-k)$.

\justify
Let $(u_\nu)_{\nu \in \mathbb{N}}$ be a bounded sequence in $H^{\ell,q}_G(M)$. 
By the density of $\mc{D}_G(M)$ in $H^{\ell,q}_G(M)$, for each $\nu \in \mathbb{N}$ we may choose a function $z_\nu \in \mc{D}_G(M)$ such that
$$
\norm{u_\nu - z_\nu}_{H^{\ell,q}(M)} < \frac{1}{\nu}.
$$
Then $(z_\nu)_{\nu \in \mathbb{N}}$ is also a bounded sequence in $H^{\ell,q}_G(M)$.
By \ref{inequacion caso entero lema embeding accion}, for each $m \in \{1,\dots,N\}$ the corresponding sequence $(\widetilde{z}_{\nu,m})_{\nu \in \mathbb{N}}$ is bounded in $H^{\ell,q}(V_m)$. By the Rellich--Kondrachov theorem on $V_m$ (extracting successively $N$ times across the finite cover), there exists a common subsequence $(z_{\nu_\alpha})_{\alpha \in \mathbb{N}} \subset (z_\nu)_{\nu \in \mathbb{N}}$ such that $(\widetilde{z}_{\nu_{\alpha},m})_{\alpha}\subset (\widetilde{z}_{\nu,m})_\nu$ is a Cauchy sequence in $L^p(V_m)$ for each $m=1,\dots,N$. Since 
\begin{equation}  \norm{z_{\nu_{\alpha}} - z_{\nu_{\beta}} }_{L^p(M)} \leq \sum_{m=1}^N K_m \norm{\widetilde{z}_{\nu_{\alpha},m} - \widetilde{z}_{\nu_{\beta},m}}_{L^p(V_m)} \xrightarrow{\alpha, \beta \to \infty} 0,
\end{equation}
it follows that $(z_{\nu_\alpha})_{\alpha \in \mathbb{N}}$ is a Cauchy sequence in $L^p(M)$ and thus converges in $L^p(M)$. 
Finally, by the continuity of the embedding $H^{\ell,q}_G(M) \hookrightarrow L^p(M)$ established above, there exists a constant $D>0$, such that the following holds
$$
\norm{u_{\nu_\alpha} - z_{\nu_\alpha}}_{L^p(M)} \leq D \norm{u_{\nu_\alpha} - z_{\nu_\alpha}}_{H^{\ell,q}(M)} < \frac{D}{\nu_\alpha} \xrightarrow{\alpha \to \infty} 0,
$$
so that the corresponding subsequence $(u_{\nu_\alpha})_{\alpha \in \mathbb{N}}$ converges to the same limit in $L^p(M)$.
\vspace{0.12cm}

We have that
\begin{equation}\label{punto 1, sobolev entero con simetria}
    H^{\ell,q}_{G}(M) \hookrightarrow L^{p}(M)
\end{equation}
is a compact embedding for $\ell \geq 1$ an integer and a real $p \in [1,\infty)$ such that $1/p>1/q-\ell/(n-k)$.

\vspace{0.1cm}

$ii)$ For the fractional Sobolev spaces $H^{s,q}_{G}(M)$ with $s,q \in \mb{R}, s>0$ non-integer, $q>1$, we seek an analogous result using the representation by complex interpolation \cite[p.25]{Taylor III}. In fact, considering $\ell \geq 0$ an integer and $\theta \in (0,1)$, one has that $H^{s, q}(M)=[H^{\ell,q}(M), H^{\ell+1,q}(M)]_{\theta}$, where $s=\ell+\theta$. Before continuing, let us consider the following claim.

\begin{claim}\label{claim:conmutatividad en espacios de interpolacion}
Let $A_{0},A_{1}$ be Banach spaces such that $A_{1} \hookrightarrow A_{0}$ continuously, and let $T: A_{0} \to A_{0}$ be a bounded linear projection whose restriction $T|_{A_{1}}: A_{1} \to A_{1}$ is also bounded, so that $T(A_{i})$ ($i=0,1$) are closed subspaces and $T(A_{1}) \hookrightarrow T(A_{0})$ continuously. Then
    \begin{equation}
[T(A_{0}),T(A_{1})]_{\theta}=T([A_{0},A_{1}]_{\theta}), \quad \theta \in (0,1).
\end{equation}
with equivalent norms.
\end{claim}

\begin{proofofclaim}
   We use $\mathcal H_{A_0,A_1}(\Omega)$ for the standard Calder\'on class on the strip \cite{Taylor I}. We prove that $T([A_{0},A_{1}]_{\theta}) \subset [T(A_{0}),T(A_{1})]_{\theta}$. Since $T: A_{i} \to A_{i}$ is linear and continuous for $i=0,1$, it follows that $T(A_{i})=\ker(I-T|_{A_{i}})$ is a closed Banach subspace of $A_{i}$ and $T: A_{i} \to T(A_{i}) \subset A_{i}$ is also continuous.
    Given $v \in [A_{0},A_{1}]_{\theta}$ there exists $u \in \mathcal{H}_{A_{0},A_{1}}(\Omega)$\,(where $\Omega = \{z \in \mathbb{C} : 0 < \operatorname{Re} z < 1\}$ denotes the open unit strip in $\mathbb{C}$) with $v=u(\theta)$. Note that $Tu: \overline{\Omega} \to T(A_{0})$ belongs 
    to $\mathcal{H}_{T(A_{0}),T(A_{1})}(\Omega)$ and satisfies $Tv=Tu(\theta)$, hence $Tv \in [T(A_{0}),T(A_{1})]_{\theta}$ with $\norm{Tv}_{[T(A_{0}),T(A_{1})]_{\theta}} \leq \max_{i=0,1}\norm{T}_{\mathcal{L}(A_{i})}\norm{v}_{[A_{0},A_{1}]_{\theta}}$.
    Conversely, given $u \in [T(A_{0}), T(A_{1})]_{\theta}$, there exists $f \in \mathcal{H}_{T(A_{0}), T(A_{1})}(\Omega)$ such that $u=f(\theta)$. Note 
    that $$\mathcal{H}_{T(A_{0}),T(A_{1})}(\Omega) \subset \mathcal{H}_{A_{0},A_{1}}(\Omega)$$
    because $T(A_{i}) \subset A_{i}$ for $i=0,1$. 
    Therefore $u \in [A_{0},A_{1}]_{\theta}$ with $\norm{u}_{[A_{0},A_{1}]_{\theta}} \leq \norm{u}_{[T(A_{0}),T(A_{1})]_{\theta}}$, and $Tu \in T([A_{0},A_{1}]_{\theta})$. 
    Using \cite[p.2.4]{Lunardi},
    we have $[T(A_{0}), T(A_{1})]_{\theta} \subset T(A_{0}) + T(A_{1}) \subset T(A_{0})$,
    then there exists $\omega \in A_{0}$ such that $u=T(\omega)$ and by the idempotence of $T$ we have $Tu=T(T(\omega))=T^{2}(\omega)=T(\omega)=u$.
\end{proofofclaim}
\justify

Since $G$ is a compact group, there exists a bi-invariant Haar measure $\sigma_G$ on $G$(see \cite[ch.9]{Cohn}).\\
Define the operator $P_{G}: H^{s,q}(M) \to H^{s,q}(M)$ for any real $s\geq 0$ and  $q > 1$ by 
\begin{equation}\label{operador de accion}
    (P_{G}u)(x) = \frac{1}{\sigma_{G}(G)}\int_{G}u(\sigma(x))d\sigma_{G}(\sigma)
\end{equation}
this operator is linear and bounded, and is also an idempotent operator.
In particular, for the endpoint $s=0$ (where $H^{0,q}(M)=L^{q}(M)$ and $H^{0,q}_{G}(M)=L^{q}_{G}(M)$), since every $\sigma \in G \leq \mathrm{Isom}_g(M)$ preserves the Riemannian volume measure $dV_g$, Minkowski's integral inequality ensures that $P_G: L^q(M) \to L^q_G(M)$ satisfies
$$
\norm{P_G u}_{L^q(M)} \leq \frac{1}{\sigma_G(G)}\int_G \norm{u \circ \sigma}_{L^q(M)}\,d\sigma_G(\sigma) = \norm{u}_{L^q(M)}.
$$
For fixed $\sigma\in G$, use right invariance of Haar measure and the change of variables $\varrho=\sigma^{\prime}\circ\sigma$, so $d\sigma_G(\varrho)=d\sigma_G(\sigma^{\prime})$. For smooth $u$, Fubini\'s theorem gives the identity below; it extends to $H^{s,q}(M)$ by density and boundedness. Then
\begin{align*}
 (P_G^2 u)(x) = P_G(P_G u)(x) &
   = \frac{1}{\sigma_{G}(G)^{2}}\int_G \left(\int_G u(\sigma^{\prime}(\sigma(x)))  d\sigma_{G}(\sigma^{\prime}) \right) d\sigma_{G}(\sigma) 
   = (P_G u)(x).
\end{align*}
Haar invariance gives $(P_Gu)\circ\tau=P_Gu$ for each $\tau\in G$, while $P_Gu=u$ for every invariant $u$. Thus the range is precisely $H^{s,q}_G(M)$, including when $s=0$. It is closed because it equals $\ker(I-P_G)$.
Note that $P_{G}(H^{s,q}(M))=H^{s,q}_{G}(M)$ and $P_{G}$ satisfies the previous Claim \ref{claim:conmutatividad en espacios de interpolacion}. Hence  we have that for any real numbers $\sigma, s \geq 0$ with $\sigma < s$, the following characterization via complex interpolation holds
\begin{equation}\label{interpretacion de sobolev con accion bajo interpolacion}
  [H^{\sigma,q}_{G}(M), H^{s,q}_{G}(M)]_{\theta}=H_{G}^{\theta s +(1-\theta)\sigma,q}(M). 
\end{equation}
In particular, by complex interpolation and the symmetric property \cite[d.2.3]{Lunardi}, we can express the fractional Sobolev space $H^{s,q}_{G}(M)$ as 
\begin{equation}
   H^{s, q}_{G}(M)=[H^{\ell,q}_{G}(M), H^{\ell+1,q}_{G}(M)]_{\theta}=[H^{\ell+1,q}_{G}(M), H^{\ell,q}_{G}(M)]_{1-\theta} 
\end{equation}
where $s=\ell+\theta$ with $\ell \in \mb{Z}_{\geq 0}$ and $\theta \in (0,1)$. 
Note that the following computation holds:
\begin{equation}\label{truco algebraico lema final tesis}
    \mfrac{1}{q}-\mfrac{s}{n-k}=\mfrac{1}{q}-\mfrac{\ell+\theta}{n-k}=(1-\theta)\left( \mfrac{1}{q}-\mfrac{\ell}{n-k}\right) + \theta\left(\mfrac{1}{q}-\mfrac{\ell+1}{n-k} \right).
\end{equation}
 If we consider the condition on $p \in [1,\infty)$ given by $\frac{1}{p} > \frac{1}{q} - \frac{s}{n-k}$, or its equivalent
\begin{equation}\label{condición para embeddings con simetria caso fraccionario}
   s > (n-k)\left(\mfrac{1}{q}-\mfrac{1}{p} \right),
\end{equation}
by \eqref{truco algebraico lema final tesis} we can choose some real numbers $p_{0},p_{1} \geq 1$ such that 
\begin{equation}
 \mfrac{1}{p_{0}} \geq \mfrac{1}{q} - \mfrac{\ell}{(n-k)}, \quad \mfrac{1}{p_{1}} > \mfrac{1}{q} - \mfrac{\ell+1}{n-k}, \quad   \mfrac{1}{p}=\mfrac{1-\theta}{p_{0}} + \mfrac{\theta}{p_{1}}. 
\end{equation}
Note that by \cite[p.325]{Taylor I} we can express
the space $L^{p}(M)$
as a complex interpolation couple by $[L^{p_{1}}(M),L^{p_{0}}(M)]_{1-\theta}=L^{p}(M)$
with $1\leq p < \infty$.
Since $q>1$, the integer Sobolev space 
$H^{\ell,q}(M)$ 
is reflexive for every $\ell \geq 0$ integer, and by the fact that $H^{\ell,q}_{G}(M)$ is a closed  Banach subspace of $H^{\ell,q}(M)$, because $H^{\ell,q}_{G}(M) = P_{G}(H^{\ell,q}(M))$, and $P_{G}(H^{\ell,q}(M)) = \text{ker}\,(I - P_{G})$. Since $I$ and $P_{G}$ are linear and continuous operators, it follows that $P_{G}(H^{\ell,q}(M))$ is the preimage of $\{0\}$ under a linear and continuous operator, and thus a closed subspace of $H^{\ell,q}(M)$; hence $H^{\ell,q}_{G}(M)$ is reflexive by \cite[p.3.20]{Brezis}.  
Using that 
\begin{equation}
    H^{\ell+1,q}_{G}(M)=[H^{\ell+2,q}_{G}(M), H^{\ell,q}_{G}(M)]_{1/2}
\end{equation}
and that the 
inclusion $H^{\ell,q}_{G}(M) \hookrightarrow L^{p_{0}}(M)$ is linear and bounded for any integer $\ell \geq 0$(note that for $\ell=0$ we have $L^q_G(M) \hookrightarrow L^{p_0}(M)$ bounded for $1 \leq p_0 \leq q$ since $M$ is compact), and compact if $\ell \in \mathbb{Z}_{\geq 1}$ and $1/p_0 > 1/q -\ell/(n-k)$, while the inclusion $H^{\ell+1,q}_{G}(M) \hookrightarrow L^{p_{1}}(M)$ is linear, bounded and compact for all $\ell \in \mb{Z}_{\geq 0}$ by \eqref{punto 1, sobolev entero con simetria}, we conclude by Theorem \eqref{Cwikel theorem for interpolation compacity1} that the embedding
\begin{equation}\label{embedding caso fraccionario parte ii}[H^{\ell+1,q}_{G}(M), H^{\ell,q}_{G}(M)]_{1-\theta} \hookrightarrow [L^{p_{1}}(M), L^{p_{0}}(M)]_{1-\theta}
\end{equation}
is compact for $0<\theta <1$. Therefore the embedding $H^{s,q}_{G}(M) \hookrightarrow L^{p}(M)$ is compact whenever $1/p>1/q-s/(n-k)$.

Let us prove that the embedding \eqref{embedding caso fraccionario parte ii} is continuous for the critical case $p_{*}=\frac{q(n-k)}{(n-k)-sq}$. For this, we use the smooth function $\widetilde{u}_{m}$ from item \hyperref[punto 1, sobolev entero con simetria]{$i)$}. For each $m\in\{1,\dots,N\}$, define the bounded linear operator $\mathcal{R}_{m}:L^{q}(M)\to L^{q}(V_{m})$ by 
$$
(\mathcal{R}_{m}u)(y)=\frac{1}{\mathrm{vol}_{\overline{g}}(U_{m})}\int_{U_{m}}(P_{G}u)\bigl(\varphi_{m}^{-1}(x,y)\bigr)\,dx.
$$
Let $u \in \mc{D}_{G}(M)$. Using the $G$-invariance, $P_{G}u=u$ and hence $\mc{R}_{m}u=\widetilde{u}_{m}$. 
By the linearity of the integral and of the operator $P_{G}$, it follows that $\mc{R}_{m}$ is linear. 
Let us prove that $\mc{R}_{m}$ is a continuous operator from $H^{\vartheta,q}(M)$ to $H^{\vartheta,q}(V_{m})$ for any integer $\vartheta \geq 0$. Indeed, using inequality \eqref{inequacion caso entero lema embeding accion}, there exists a constant $C_{m,\vartheta,q}>0$ such that
\begin{align*}
    \norm{\mc{R}_{m}u}^{q}_{H^{\vartheta,q}(V_{m})} &= \sum_{0 \leq |j| \leq \vartheta}\int_{V_{m}}|\nabla^{j}\mc{R}_{m}u(y)|^{q}\,dy = \sum_{0 \leq |j| \leq \vartheta}\int_{V_{m}}|\nabla^{j}\widetilde{u}_{m}(y)|^{q} \, dy \leq C_{m,\vartheta,q}\norm{u}^{q}_{H^{\vartheta,q}(M)}.
\end{align*}
In particular, $\mc{R}_{m}$ is a continuous operator for $\vartheta \in \{\ell, \ell +1\}$. Since $V_{m} \subset \mb{R}^{n-k_{m}}$ is a bounded open domain with smooth boundary, by Theorem \cite[p.2.1]{Taylor I}, complex interpolation between $j=\ell$ and $j=\ell+1$ yields that there exists $C_{m,s,q}>0$ such that
\begin{equation}\label{estimacion perfil fraccionario}
\norm{\widetilde{u}_{m}}_{H^{s,q}(V_{m})}=\norm{\mathcal{R}_{m}u}_{H^{s,q}(V_{m})}\leq C_{m,s,q}\norm{u}_{H^{s,q}(M)}.
\end{equation}
Define $s_{m}=s(n-k_{m})/(n-k)$; then $0<s_{m}\leq s$, $s_{m}q < d_{m}$ and $p_{*}=\frac{q(n-k_{m})}{(n-k_{m})-s_{m}q}$. Therefore, using that $V_{m}$ is a smooth bounded extension domain, by Theorem \cite[ch.13,p.6.4]{Taylor I}, and estimate \eqref{estimacion perfil fraccionario}, we obtain the following relation
\begin{align}\label{relacion de embeddings caso ii critico}
    H^{s,q}(M) \xrightarrow{\mc{R}_{m}} H^{s,q}(V_{m}) \hookrightarrow H^{s_{m},q}(V_{m}) \hookrightarrow L^{p_{*}}(V_{m}). 
\end{align}
Taking into account a smooth non-negative partition of unity $(\varrho_{m})_{m=1}^{N}$ subordinate to the cover $\{(\Omega_{m}, \varphi_{m})\}_{m=1}^{N}$ and inequality \eqref{estimacion perfil fraccionario}, we have
\begin{align}
\norm{u}_{L^{p_{*}}(M)}=\left( \sum_{m=1}^{N} \int_{\Omega_{m}}\varrho_{m}|u|^{p}\,dV_{g}  \right)^{\frac{1}{p}}&\leq\sum_{m=1}^{N}K_{m}\norm{\widetilde{u}_{m}}_{L^{p_{*}}(V_{m})}\leq C\norm{u}_{H^{s,q}(M)}.
\end{align}
for some constant $C>0$. By the density of $\mathcal{D}_{G}(M)$ in $H^{s,q}_{G}(M)$ (Lemma \ref{espacio denso en espacios de sobolev fraccionarios}), this proves the continuous embedding $H^{s,q}_{G}(M)\hookrightarrow L^{p_{*}}(M)$.
\vspace{.12cm}

$iii)$ Assume $sq>n-k$ and fix $0<\eta<\min\{1,s-\frac{n-k}{q}\}$.
Proceeding similarly to the proof of item \hyperref[punto 1, sobolev entero con simetria]{$i)$}. Let $u \in \mc{D}_{G}(M)$. 
Since $(M,g)$ is compact, there exists a finite open covering $\mc{A}=\{(\Omega_{m}'', \varphi_{m})\}_{m=1}^{N}$ satisfying the conditions of Lemma \eqref{lemma de las cartas} such that $\varphi_{m}(\Omega_{m}'')=U''_{m} \times V''_{m}$, where $U''_{m} \subset \mb{R}^{k_{m}}$ and $V_{m}'' \subset \mb{R}^{n-k_{m}}$ are bounded open subsets satisfy 
$$U''_{m} \Subset U_{m}' \Subset U_{m}, \quad \text{and} \quad V''_{m} \Subset V_{m}' \Subset V_{m},$$
with $U_{m}$ and $V_{m}$ being the bounded open subsets defined in \eqref{punto 1, sobolev entero con simetria 1}, and $U_{m}', V_{m}'$ bounded open sets satisfying the same respective properties as $U_{m}$ and $V_{m}$, inducing a finite open cover of $M$ given by $\mc{B}=\{(\Omega_{m}', \varphi_{m})\}_{m=1}^{N}$ and $\varphi_{m}^{-1}(U_{m}' \times V_{m}')=\Omega_{m}'$.
We know from \hyperref[punto 1, sobolev entero con simetria]{$i)$}, that there exists for every $x,x'\in U_{m}$ and $y \in V_{m}$, a smooth $G$-invariant function $\widetilde{u}_{m} \in C^{\infty}(\mathbb{R}^{n-k_{m}})$, defined by 
\begin{equation}
  \widetilde{u}_{m}(y)=u \circ \varphi_{m}^{-1}(x,y). 
\end{equation}
Let $\{\psi_{m}\}_{m=1}^{N}$ be a family of smooth cut-off functions subordinate to $\mc{B}$, whose existence is justified in \cite[p.2.25]{Lee}
such that $\text{supp}\,(\psi_{m}) \Subset \Omega_{m}'$.\\
For a fixed $m$, let $\rho_{m} \in C^{\infty}_{0}(U_{m}')$ satisfy $\rho_{m}\geq 0$ and $\int_{U_{m}'}\rho_{m}(x)\, dx=1$, such that $\psi_{m} \equiv 1$ on a neighborhood of  $\varphi_{m}^{-1}(\text{supp}\,(\rho_{m}) \times \overline{V_m''})$ containing $\overline{\Omega_{m}''}$. 
Define the operator $T_{m}$ as 
\begin{equation}
    T_{m} : C^{\infty}_{0}(M) \to C^{\infty}_{0}(V_{m}'), \quad u \mapsto T_{m}u(y)=\int_{U_{m}'}\rho_{m}(x)((\psi_{m}u) \circ \varphi_{m}^{-1})(x,y) \,dx.
\end{equation}
Note that on $V_m''$,  
$T_{m}u(y) =\int_{\text{supp}\,(\rho_{m})}\rho_{m}(x)\widetilde{u}_{m}(y)dx=\widetilde{u}_{m}(y).$
\vspace{0.1cm}

Now we prove that $T_{m}$ is bounded on $H^{s,q}(M)$. To this end, we first show boundedness on integer-order Sobolev spaces $H^{j,q}(M)$ for $j \geq 0$. By Hölder's inequality with $q'$ conjugate of $q$, there exist $C_{m}, \widetilde{C}_{m}>0$ constants such that
\begin{align*}
\norm{T_{m}u}_{H^{j,q}(V_{m}')}^{q}
&\leq  \left(\sum_{|\alpha|\leq \,j}\int_{V_{m}'} \left( \int_{U'_{m}}|\rho_{m}(x)||\partial_{y}^{\alpha}(\psi_{m}u \circ \varphi_{m}^{-1})(x,y)|\,dx \right)^{q} \,dy    \right) \\
&\leq  \left(\sum_{|\alpha|\leq \,j}\norm{\rho_{m}}_{L^{q'}(U_{m}')}^{q} 
\int_{U_{m}' \times V_{m}'}|\partial_{y}^{\alpha}(\psi_{m}u \circ \varphi_{m}^{-1})(x,y)|^{q}\, dx dy    \right) \\
&\leq C_{m}\sum_{|\alpha|\leq \,j}\norm{\partial_{y}^{\alpha}(\psi_{m}u \circ \varphi_{m}^{-1})  }^{q}_{L^{q}(U_{m}' \times V_{m}')} \leq \left(\frac{C_{m}}{\widetilde{B}_{m,0}}\right)\norm{u}^{q}_{H^{j,q}(\Omega_{m}')}\leq \widetilde{C}_{m}\norm{u}_{H^{j,q}(M)}^{q}.
\end{align*}
Here  $\widetilde{B}_{m,0}>0$ is a constant arising from the inequality \eqref{inequacion caso entero lema embeding accion} satisfying 
$$\norm{\partial_{y}^{\alpha}(\psi_{m}u \circ \varphi_{m}^{-1})}^{q}_{L^{q}(U_{m}' \times V_{m}')} \leq \frac{1}{\widetilde{B}_{m,0}}\int_{\Omega_{m}'}|\partial^{\alpha}_{y}u|^{q}dV_{g}.$$
Thus, $T_{m}$ is bounded on integer-order Sobolev spaces $j \geq 0$ for $q>1$, so
$$T_{m}: H^{j,q}(M) \to H^{j,q}(V_{m}')$$
is linear and continuous. 
Since $V_{m}'$ is a bounded open set with smooth boundary, fractional Sobolev spaces on $V_{m}'$ and on the compact manifold $M$ can be defined via complex interpolation. By Theorem \cite[p.2.1;323]{Taylor I}, for $\theta \in (0,1)$, integer $\ell\geq 0$, and $s=\ell+\theta$, the operator
\begin{equation}
    T_{m}:  H^{s,q}(M):=[H^{\ell,q}(M), H^{\ell+1,q}(M)]_{\theta} \to H^{s,q}(V_{m}'):=[H^{\ell,q}(V_{m}'), H^{\ell+1,q}(V_{m}')]_{\theta}
\end{equation}
is linear and continuous. Consequently,
\begin{equation}\label{desigualdad lema principal 1} \norm{T_{m}u}_{H^{s,q}(V_{m}')}\leq \norm{u}_{H^{s,q}(\Omega^{'}_{m})} \leq \norm{u}_{H^{s,q}(M)}.
\end{equation}
Since $V'_{m}$ is a smooth bounded extension domain, the Sobolev--Morrey embedding \cite[p.8.4]{Taylor III} gives $H^{s,q}(V'_{m}) \hookrightarrow C^{0,\eta}(\overline{V'_{m}})$ for $s>\eta + (n-k_{m})/q$ and $0<\eta < 1$. This follows by considering the extension operator $\mc{E}_{m}$ on $\mb{R}^{n-k_{m}}$ and the restriction operator $\mc{R}_{m}$ associated with $V'_{m}$; indeed, we have the following chain
\begin{align*}
    H^{s,q}(V'_{m}) \xrightarrow{\mc{E}_{m}} H^{s,q}(\mb{R}^{n-k_{m}}) \hookrightarrow H^{\eta + \frac{n-k_{m}}{q}, q}(\mb{R}^{n-k_{m}}) \hookrightarrow C^{0,\eta}(\mb{R}^{n-k_{m}}) \xrightarrow{\mc{R}_{m}} C^{0,\eta}(\,\overline{V'_{m}}).
\end{align*}
Therefore the initial choice $0<\eta < \min\{1, s- \frac{n-k}{q}\}$ satisfies this condition for all $m$. This applies for every real $s>0$, including when $s$ is an integer. Hence there exists $K>0$ such that 
\begin{equation}\label{desigualdad lema principal 2}
\norm{\widetilde{u}_{m}}_{C^{0,\eta}(V_{m}'')} \leq \norm{T_m u}_{C^{0,\eta}(V_m')} \leq K^{1/q} \norm{T_{m}u}_{H^{s,q}(V_{m}')}.  
\end{equation}
We now prove that for some constant $C_{m,\eta}>0$; we have that
\begin{equation}\label{desigualdad lema principal 3}
\norm{u}_{C^{0,\eta}(\,\overline{\Omega_{m}''}\,)} \leq C_{m,\eta}\norm{u}_{H^{s,q}(M)}.
\end{equation}
Let $p,p' \in \overline{\Omega''_{m}}$, then $p=\varphi_{m}^{-1}(x,y)$,  $p'=\varphi_{m}^{-1}(x',y')$ for some  $(x,y), (x',y') \in \overline{U_{m}'' \times V_{m}''}$. 
Since $\Omega''_{m} \Subset \Omega'_{m}$, $\overline{\Omega''_{m}}$ is a compact subset of $\Omega'_{m}$. By  Theorem \cite[t.13.29]{Lee}, there exist constants $c_{m},d_{m}>0$ such that for any piecewise smooth curve $\gamma$ lying entirely in $\overline{\Omega''_{m}}$ we have
\begin{equation*}
c_{m}\ell_{\overline{g}}(\gamma) \leq \ell_{g}(\gamma) \leq d_{m}\ell_{\overline{g}}(\gamma).
\end{equation*}
Let us first consider curves whose image lies in $\Omega'_{m}$. Let $\gamma_{c}:(-\varepsilon, \varepsilon) \to \Omega'_{m}$ with $\varepsilon>1$ be a piecewise smooth curve connecting $p$ and $p'$ with $\gamma_{c}(0)=p$ and $\gamma_{c}(1)=p'$, and let $\widetilde{\gamma_{c}}=\varphi_{m} \circ \gamma_{c}$ be its coordinate representation in $(\Omega''_{m},\varphi_{m})$. In fact
\begin{align*}
    \ell_{g}(\gamma_{c})=\int_{0}^{1}|\gamma_{c}'(t)|_{g}\,dt \geq c_{m}\int_{0}^{1}|\widetilde{\gamma_{c}}'(t)|_{\overline{g}}\, dt \geq c_{m}|\widetilde{\gamma_{c}}(1) - \widetilde{\gamma_{c}}(0)| \geq c_{m}|y-y'|.
\end{align*}
Taking the infimum over all curves $\gamma_{c}$ in $\Omega'_{m}$, we have $ d_g(p,p') = \inf_{\gamma_c} \ell_g(\gamma_c) \geq c_m|y-y'|$. Note that this is valid only for curves $\gamma_{c}$ that lie inside $\Omega_{m}'$. 
Hence $|y-y'|\leq c_{m}^{-1}d_{g}(p,p')$. To pass to the global Riemannian distance, we also consider curves that leave this chart.\\
In the case where we have a curve $\gamma_{c}$ that leaves $\Omega'_{m}$, we use the fact that $\Omega''_{m} \Subset \Omega'_{m}$, hence $$r_m := d_g\bigl(\,\overline{\Omega_m''}, M \setminus \Omega_m'\bigr) > 0.$$ 
This implies that $\ell_{g}(\gamma_{c})\geq r_{m}$; moreover, since $y,y' \in \overline{V''_{m}}$ and $V''_{m}$ is a bounded subset of $\mb{R}^{n-k_{m}}$, we have $|y-y'|\leq \mathrm{diam}(V''_{m})<\infty$. 
Therefore, for this case we have the estimate $\ell_g(\gamma_c) \geq r_m \geq \frac{r_m}{\mathrm{diam}(V''_{m})}|y-y'|.$
In fact, taking $\widetilde{K}_{m}=\max\left\{\frac{1}{c_{m}}, \frac{\mathrm{diam}(V''_{m})}{r_{m}}   \right\}$, we obtain the following estimate by taking the infimum over all piecewise $C^1$ connecting curves, using the two cases above: 
$$|y-y'| \leq \widetilde{K}_{m}d_{g}(p,p').$$

\justify

By inequality \eqref{desigualdad lema principal 2}, we have that $\norm{\widetilde{u}_{m}}_{C^{0,\eta}(\,\overline{V_{m}''}\,)} \leq \|T_{m}u\|_{C^{0,\eta}(V_{m}')}$, and thus

\justify

\begin{align*}
[\widetilde{u}_{m}]_{C^{0,\eta}(V''_{m})}=\sup_{\substack{y,y' \in V_{m}'' \\ y \neq y'}}\left\{\frac{|\widetilde{u}_{m}(y)-\widetilde{u}_{m}(y')|}{|y-y'|^{\eta}} \right\} &\leq \norm{u}_{H^{s,q}(M)}.  
\end{align*}
Therefore, there exists a constant $C_{m,\eta}>0$ such that
\begin{equation}\label{pre desigualdad lema principal 3.1}
    |u(p) - u(p')|=|\widetilde{u}_{m}(y)-\widetilde{u}_{m}(y')| \leq C_{m,\eta}\norm{u}_{H^{s,q}(M)}d_{g}(p,p')^{\eta}
\end{equation}
implying that 
\begin{align}\label{desigualdad lema principal 3.1}
  \norm{u}_{C^{0,\eta}(\Omega_{m}'')} \leq \norm{u}_{C^{0,\eta}(\,\overline{\Omega_{m}''}\,)} \leq C_{m,\eta}\norm{u}_{H^{s,q}(M)}. 
\end{align}
Since $\mc{A}$ is an open cover of $M$, inequality \eqref{desigualdad lema principal 3.1} yields
\begin{equation*}
\norm{u}_{L^{\infty}(M)}=\max_{1\leq m \leq N}\left\{\norm{u}_{L^{\infty}(\Omega''_{m})} \right\} \leq \left( \max_{1 \leq m \leq N} C_{m,\eta}\right)\norm{u}_{H^{s,q}(M)}.
\end{equation*}
Now we extend the inequality \eqref{desigualdad lema principal 3.1} globally to show that
$\norm{u}_{C^{0,\eta}(M)} \leq C_{\eta}\norm{u}_{H^{s,q}(M)}$, for some constant $C_{\eta}>0$. 
The metric space $(M,d_{g})$ induced by the metric $g$ is sequentially compact, so it admits a Lebesgue number for the cover $\mc{A}$. Let $\delta>0$ be the Lebesgue number of $\mc{A}$, so that for any $p \in M$, there exists some $m_{0} \in \{1,...,N\}$, such that $B_{g}(p,\delta) \subseteq \Omega_{m_{0}}''$, with $B_{g}(p,\delta)$ the open geodesic ball centered at $p$ with radius $\delta$.
For any $p,p' \in M$ with $p\neq p'$, we consider two cases:

\vspace{0.12cm}

$a')$ $d_{g}(p,p') < \delta$: Then $p,p' \in \Omega''_{m_{i}}$ for some $m_{i}$, and by inequality \eqref{pre desigualdad lema principal 3.1},
\begin{align*}
    \frac{|u(p)-u(p')|}{d_{g}(p,p')^{\eta}} &\leq C_{m_{i}, \eta}\norm{u}_{H^{s,q}(M)} \leq  \left( \max_{1 \leq m \leq N} C_{m,\eta}\right) \norm{u}_{H^{s,q}(M)}.
\end{align*}

\vspace{0.12cm}

$b')$ $d_{g}(p,p') \geq \delta$: In this case, $p$ and $p'$ may lie in different charts. Using the $L^{\infty}(M)$ norm,
$$|u(p) - u(p')| \leq |u(p)| + |u(p')| \leq 2 \norm{u}_{L^{\infty}(M)}$$
Setting $C_{\infty}=\max_{1 \leq m \leq N} C_{m,\eta}$, we have that
\begin{align*}
    \frac{|u(p)-u(p')|}{d_{g}(p,p')^{\eta}} &\leq \frac{2\norm{u}_{L^{\infty}(M)}}{\delta^{\eta}} \leq \left(\frac{2C_{\infty}}{\delta^{\eta}}    \right)\norm{u}_{H^{s,q}(M)}.
\end{align*}
Hence, 
\begin{equation}
    [u]_{C^{0,\eta}(M)} \leq \left(\max\left\{\frac{2C_{\infty}}{\delta^{\eta}}, \left(\max_{1\leq m \leq N}C_{m,\eta}\right)\right\} \right)\norm{u}_{H^{s,q}(M)}. 
\end{equation}
Thus, there exists a constant $C_{\eta}>0$ such that 
\begin{equation}\label{desigualdad lema principal 3.2}
\norm{u}_{C^{0,\eta}(M)} \leq C_{\eta}\norm{u}_{H^{s,q}(M)}.
\end{equation}
Since $\mc{D}_{G}(M)$ is dense in $H^{s,q}_{G}(M)$ by Lemma \eqref{espacio denso en espacios de sobolev fraccionarios}, equation \eqref{desigualdad lema principal 3.2} holds for all $u \in H^{s,q}_{G}(M)$, yielding the continuous embedding
\begin{equation}\label{embedding lema final caso morrey}
   H^{s,q}_{G}(M) \hookrightarrow C^{0,\eta}(M) 
\end{equation}
for $0< \eta < \min\left\{1, s-\frac{n-k}{q}\right\}$.

\vspace{0.12cm}
Let $(u_{n}) \subset H^{s,q}_{G}(M)$ be a bounded sequence. Then $(u_{n}) \subset C^{0,\eta}(M)$ is a bounded and continuous sequence of functions. Then there exists a family $\mc{F} \subseteq \,C^{0}((M,d_{g}),\mb{R})$ which is equicontinuous and pointwise bounded, such that $(u_{n}) \subset \mc{F}$. 
Then by the Arzelà-Ascoli theorem, there exists $(u_{n_{k}}) \subset (u_{n})$ that converges uniformly on $M$ in the supremum norm $\|\cdot\|_{\infty}$ to some continuous function $u \in C^0(M)$.
Therefore the embedding $(C^{0,\eta}(M),\|\cdot\|_{C^{0,\eta}(M)}) \hookrightarrow (C^{0}(M), \|\cdot\|_{\infty})$ is compact, which implies that
\begin{equation}\label{embedding lema final caso morrey 1}
   H^{s,q}_{G}(M) \hookrightarrow C^{0}(M) 
\end{equation}
is a compact embedding. 
Since the embedding $C^{0}(M) \hookrightarrow L^{\infty}(M)$ is continuous, it follows that
\begin{equation}\label{embedding lema final caso morrey 2}
    H^{s,q}_{G}(M) \hookrightarrow L^{\infty}(M)
\end{equation}
is a compact embedding for $sq>n-k$.
\end{proof}

\vspace{0.12cm}
\justify

Let $f$ be a nonnegative real-valued function that satisfies the ellipticity condition (\ref{codicion de elipticidad de una fucnion real valuada para Riemannian manifolds}). We introduce the space $\widetilde{H}_{f,G}(M,g)$ for $G \leq \text{Isom}_{g}(M)$ a closed Lie subgroup
\begin{equation}\label{espacio de Sobolev en verson sumas con acción}
    \widetilde{H}_{f, G}(M,g) = \left\{ h \in \widetilde{H}_{f}(M,g) \mid h \circ \sigma=h, \forall \sigma \in G \right\}.
\end{equation}
Similarly we consider the Banach space $L^{p}_{G}(M)$ as $$L^{p}_{G}(M)=\left\{u \in L^{p}(M) \mid u \circ \sigma=u, \forall \sigma \in G \right\}.$$

\begin{lemma}\label{sistema lineal caso riemanniano con acción}
    For every $h \in L^2_{G}(M,g)$, the linear problem
\begin{equation}\label{ecuacion lineal caso Riemanniano}
f(\Delta_g)\phi + \phi = h
\end{equation}
has a unique solution $\phi \in \widetilde{H}_{f,G}(M,g)$. Moreover,
$\norm{\phi}_{\widetilde{H}_{f}(M,g)}=\norm{h}_{L^{2}(M,g)}$.
\end{lemma}

\begin{proof}
Note that $\Delta_{g}$ commutes with the action of $G$, i.e., for any $\psi \in G$ and any function $f$ on $M$ one has
$\Delta_{g}(f \circ \psi)=(\Delta_{g}f)\circ \psi$, therefore
$f(\Delta_{g})$ inherits commutation with $G$, that is, $(f(\Delta_{g})\phi)\circ 
\gamma=f(\Delta_{g})(\phi \circ \gamma)$ for every $\gamma \in G$. Hence if $\phi \in 
\widetilde{H}_{f,G}(M,g)$ then $h=f(\Delta_{g})\phi + \phi$ is 
also invariant. If $h \in L^{2}(M,g)$ is invariant with respect to $G$, the unique solution $\phi$ of (\ref{ecuacion lineal caso Riemanniano}) given by (\ref{unica solucion caso lineal 
caso Riemanniano}) is also invariant, because, if $h$ is $G$-invariant, the operator $L_{c}$ is also $G$-invariant, then $L_{c}(\phi \circ \gamma)=(L_{c}\phi) \circ \gamma=h \circ \gamma=h$ for every $\gamma \in G$, then  by uniqueness of the solution, $\phi \circ \gamma=\phi$. The equality of the norms follows from (\ref{sistema lineal caso riemanniano}) and is given by
\begin{equation*}
\norm{\phi}^{2}_{\widetilde{H}_{f}(M,g)}= \sum |\phi_{i}|^{2}(1+f(\lambda_{i}))^{2}=\sum |h_{i}|^{2}=\sum_{i=1}^{\infty}\sum_{j=1}^{\infty}h_{i}h_{j}\delta_{ij}=    \norm{h}^{2}_{L^{2}(M,g)}
\end{equation*}
\end{proof}

Having established the compact Sobolev embeddings under group actions in Theorem \ref{Lema principal caso con accion} together with the well-posedness of the $G$-invariant linear problem in Lemma \ref{sistema lineal caso riemanniano con acción}, we now have all the necessary tools to address the nonlinear problem. We are thus in a position to state and prove the main existence result of this section.

\newpage

\begin{theorem}\label{Teorema Final}
    Let $(M,g)$ be a compact, connected and oriented Riemannian $n$-manifold without boundary. Let $G$ be a closed subgroup of $\mathrm{Isom}_{g}(M)$ whose action on $M$ is non-transitive. 
    Suppose that $|\mathcal{O}_{G}(x)|=+\infty, \forall x \in M$, and $k=\min_{x \in M}\mathrm{dim}~\mathcal{O}_{G}(x)$. Assume that $f$ satisfies (\ref{codicion de elipticidad de una fucnion real valuada para Riemannian manifolds}),
    $U: M \times \mathbb{R} \to \mathbb{R}$ is a Carathéodory non-linearity which is $G$-invariant with $U(\cdot, 0) \not\equiv 0$ and
    $C^{1}$ in the second variable, 
    and there exist constants $\alpha>1$, $\beta \geq 1$, $C > 0$, and a function $h \in L_{G}^2(M,g)$ such that the following holds:
\begin{equation}\label{desigualdad de crecimiento}
|U(x, y) + y| \leq C \left( |h(x)| + |y|^\alpha \right),\quad \left| \frac{\partial}{\partial y} \left( U(x, y) + y \right) \right| \leq C ( 1 + |y|^\beta ).
\end{equation}
Then, there exists $\varepsilon \in (0,1]$ and $0<\rho_{\varepsilon} < 1$ such that whenever $\norm{h}_{L^2(M)} < \rho_{\varepsilon}$, there exists a solution $\phi \in \widetilde{H}_{f,G}(M,g)$ to the non-linear problem \begin{equation}\label{ecuacion caso general fraccionario con accion}
f(\Delta_g)\phi - U(x, \phi) = 0.
\end{equation}
\end{theorem}

\begin{proof} 
Using the operator $L_{c}=f(\Delta_{g}) + \mathrm{Id}$ defined in (\ref{operador caso Riemanniano}) we have that solving the nonlinear equation 
(\ref{ecuacion caso general fraccionario con accion}) is equivalent to solving
\begin{equation}\label{sistema lineal caso Riemaniano con accion}
    L_{c}u=\mathcal{A}(\cdot, u)
\end{equation}
where $\mathcal{A}(\cdot,u)=U(\cdot, u) + u$. Let $s, s_{\alpha} \in \mathbb{R}_{>0}$ be such that $ \frac{n-k}{2}\left( \frac{\alpha -1}{\alpha}  \right) < s_{\alpha} < \frac{n-k}{2} <s$\label{desigualdades convenientes caso con accion}
holds; this inequality is not empty because $k<n$ because the action $G \curvearrowright M$ is non-transitive. 
Note that for this real $s>0$ such that by \eqref{embedding caso riemaniano} the embedding $\widetilde{H}_{f}(M,g) \hookrightarrow H^{s}(M)$ holds, which implies that $\widetilde{H}_{f,G}(M,g) \subset \widetilde{H}_{f}(M,g) \subset H^{s}(M)$.
Let us prove some results that will serve as tools later on.
The first of these is the following. Let $t>0$ be a real number, let us prove that 
\begin{equation}\label{teorema final tesis embedding 1}
\widetilde{H}_{f,G}(M,g) \hookrightarrow H^{t}_{G}(M).
\end{equation}
Since $G \leq \mathrm{Isom}_{g}(M)$ is closed and $M$ is compact, it follows that $\mathrm{Isom}_{g}(M)$ is a compact Lie group \cite[p.249]{Hebey 2}, hence $G$ is a compact Lie group; therefore, there exists a bi-invariant Haar measure on $G$, denoted by $\sigma_{G}$ \cite[ch.9]{Cohn}. Let us define the projection operator $P_{t,G}: H^{t}(M) \to H^{t}(M)$ by
\begin{equation}
    P_{t,G}u(x) = \frac{1}{\sigma_{G}(G)}\int_{G}u(\sigma(x))d\sigma_{G}(\sigma).
\end{equation}
Notice that this operator satisfies $\mathrm{ran}(P_{t,G})=H^{t}_{G}(M)$ and is a linear and bounded projection. Hence it is a closed subspace of $H^{t}(M)$, and thus a Hilbert space.
If we consider a function $u \in H^{t}(M)$ which is $G$-invariant, we have $P_{t,G}u=u$. Therefore, we obtain the following relation
\begin{equation}
\widetilde{H}_{f,G}(M,g) \xhookrightarrow{\iota_{t}} H^{t}(M) \xrightarrow{P_{t,G}} H^{t}_{G}(M). 
\end{equation}
If we define $\psi:= P_{t,G} \circ \iota_{t}: \widetilde{H}_{f,G}(M,g) \hookrightarrow H^{t}_{G}(M)$, then $\psi$ is a continuous linear operator satisfying $\psi=\iota$. 
Therefore, the embedding \eqref{teorema final tesis embedding 1} holds and is continuous.
Since $2s>n-k$, where $1\leq k < n$. By Theorem \ref{Lema principal caso con accion} for $1 \leq p < \infty$, the chain of embeddings
\begin{equation}\label{teorema final tesis embedding 2}
    \widetilde{H}_{f,G}(M,g) \hookrightarrow H^{s}_{G}(M) \hookrightarrow L^{p}(M)
\end{equation}
is continuous.
A second result is the following. 
Let $v,w\geq 0$ be real numbers with $v>w$, let us prove that
\begin{equation}\label{teorema final tesis embedding 3}
    H^{v}_{G}(M) \hookrightarrow H^{w}_{G}(M)
\end{equation}
is a compact embedding. Since $M$ is compact, using the fact that $v-\frac{n}{2}>w-\frac{n}{2}$, by Rellich-Kondrachov's Theorem \cite[t.2.9]{Hebey 2}, it follows that $H^{v}(M) \hookrightarrow H^{w}(M)$ is a compact embedding, which implies that $\iota\vert{}_{H^{v}_{G}}: H^{v}_{G}(M) \hookrightarrow H^{w}(M)$ is also compact.
Considering the bounded linear projection operator $P_{w,G}$, if we define the operator $\vartheta:=P_{w,G} \circ \iota\vert{}_{H^{v}_{G}}: H^{v}_{G}(M) \to H^{w}_{G}(M)$, it turns out to be a compact operator since it is the composition of a continuous linear operator $P_{w,G}$ with the compact one $\iota\vert{}_{H^{v}_{G}}$. Therefore, we have the following relation of embeddings
\begin{equation*}
H^{v}_{G}(M)\xhookrightarrow{\iota|_{H^{v}_{G}}} H^{w}(M) \xrightarrow{P_{w,G}} H^{w}_{G}(M)
\end{equation*}
which implies that $\vartheta=\iota$; indeed, the embedding \eqref{teorema final tesis embedding 3} is compact.
Considering $s>0$ real, there exists a real $r>0$ with $r>s$ and $s=\theta r$ for $\theta \in (0,1)$. Therefore, by the embedding \eqref{teorema final tesis embedding 3} and the fact that if $M$ is compact, then $vol_{g}(M)<\infty$, which implies that $L^{p}(M) \subset L^{q}(M)$ for all real $p,q>0$ with $p>q$, it follows that
\begin{equation}\label{teorema final tesis embedding 4}
H^{r}_{G}(M) \hookrightarrow H^{s}_{G}(M) \hookrightarrow H^{s_{\alpha}}_{G}(M) \hookrightarrow L^{2\alpha}(M) \hookrightarrow L^{2}(M)
\end{equation}
\begin{equation}\label{teorema final tesis embedding 5}
\widetilde{H}_{f,G}(M,g)
\hookrightarrow H^{r}_{G}(M) \hookrightarrow H^{s}_{G}(M)
\end{equation}
are continuous embeddings.
Note that $H^{s_{\alpha}}_{G}(M) \hookrightarrow L^{2\alpha}(M)$ is a compact embedding by Theorem \ref{Lema principal caso con accion}, since $s_{\alpha}<\frac{n-k}{2}$ and $1\leq 2\alpha=\tfrac{2(n-k)}{(n-k) - 2\left(\frac{n-k}{2}\left(\frac{\alpha-1}{\alpha} \right) \right)} < \tfrac{2(n-k)}{(n-k) - 2s_{\alpha}}.$
\justify
By the same procedure to obtain the chain \eqref{teorema final tesis embedding 4}, we can obtain a chain that will be useful below. We define $\tau=\max\{\alpha, \beta +1\}$, and a real $s_{\tau} > 0$ such that $\frac{n-k}{2}\left( \frac{\tau -1}{\tau} \right) < s_{\tau} < \frac{n-k}{2} < s$, then we have the following chain of embeddings:
\begin{equation}\label{teorema final tesis embedding 6}
       H^{s}_{G}(M) \hookrightarrow H^{s_{\tau}}_{G}(M) \hookrightarrow L^{2\tau}(M) \hookrightarrow L^{2(\beta + 1)}(M) \hookrightarrow L^{2\beta}(M)
\end{equation}
is a continuous embedding chain. 
Taking $u \in H^{s}_{G}(M)$, we want to prove that $\mathcal{A}
(\cdot,u) \in L^{2}_{G}(M,g)$. Since $u$ is $G$-invariant and $U(\cdot,u)$ is also $G$-invariant, it follows that $\mathcal{A}(\cdot,u)$ is $G$-
invariant.
We need to prove that $\mathcal{A}(\cdot,u) \in L^{2}(M)$, to get this we use (\ref{desigualdad de crecimiento}), the chain \eqref{teorema final tesis embedding 4}, the Young's inequality for the terms 
$|h|, |u|^{\alpha}$ and $\|u\|_{L^{2\alpha}(M)}^{\alpha}=\||u|^{\alpha
}\|_{L^{2}(M)}$. In fact we have that
\begin{align*}
    \norm{\mathcal{A}(\cdot,u)}_{L^{2}(M)} &\leq \sqrt{C}\left( \int_{M} (|h(x)| + |u(x)|^{\alpha})^{2}   dV_{g}\right)^{\frac{1}{2}}\leq \widetilde{C}\left(\int_{M}(|h(x)|^{2}+|u(x)|^{2\alpha})dV_{g}\right)^{\frac{1}{2}}\\
    &\leq \widetilde{C}\left(\norm{h}_{L^{2}(M)} +\norm{u}^{\alpha}_{L^{2\alpha}(M)}   \right) \leq C_{s}\left(\norm{h}_{L^{2}(M)} +\norm{u}_{H^{s}(M)}^{\alpha}   \right)
\end{align*}
for suitable constants $\widetilde{C}, C_{s}>0$. Therefore $\mathcal{A}(\cdot,u) \in L^{2}_{G}(M,g)$. Define 
\begin{equation}
\mathcal{O}_{\varepsilon}=\{u \in H^{s}_{G}(M) \mid  \norm{u}_{H^{s}(M)} \leq \varepsilon\}
\end{equation}
for an arbitrary $\varepsilon>0$. 
Note that $\mathcal{O}_{\varepsilon}$ is a closed, bounded, convex, and nonempty subset of the Hilbert space $H^{s}_{G}(M)$.
Define the map
\begin{equation}\label{operador caso con simetria}
    \mathcal{S}: \mathcal{O}_{\varepsilon} \to H^{s}_{G}(M), \quad \mathcal{S}(u)=\widetilde{u}
\end{equation}
where $\widetilde{u} \in \widetilde{H}_{f,G}(M,g)$ is the unique solution of the non-homogeneous problem $L_{c}\widetilde{u}=\mathcal{A}(\cdot,u)$. This map $\mc{S}$ is well 
defined by Lemma \ref{sistema lineal caso riemanniano con acción}.
We show that there 
exists $\varepsilon>0$ such that $\mathcal{S}: \mathcal{O}_{\varepsilon} \to \mathcal{O}_{\varepsilon}$.
Let $u \in \mathcal{O}_{\varepsilon}$.
Using the embedding  \eqref{teorema final tesis embedding 5}, Lemma \ref{sistema lineal caso riemanniano con acción} and the  values
\begin{equation*}
\varepsilon:=\min\left\{1, \left(\frac{1}{2C_{s}}   \right)^{\frac{1}{\alpha-1}}\right\} \in(0,1], \quad \rho_{\varepsilon}:=\frac{\varepsilon}{2C_{s}} \in (0,1)
\end{equation*}
where $\|h\|_{L^{2}(M)} < \rho_{\varepsilon}$. We have the following computations
\begin{align*}
    \norm{\mathcal{S}(u)}_{H^{s}(M)} &\leq \norm{\mathcal{S}(u)}_{H^{r}(M)} = \norm{\widetilde{u}}_{H^{r}(M)} \leq \norm{\widetilde{u}}_{\widetilde{H}_{f}(M)}= \norm{\mathcal{A}(\cdot, u)}_{L^{2}(M)} \\
    &\leq C_{s}\left(\norm{h}_{L^{2}(M)} +\norm{u}_{H^{s}(M)}^{\alpha}   \right)\leq C_{s}\left(\norm{h}_{L^{2}(M)} +\varepsilon^{\alpha}   \right) < \tfrac{\varepsilon}{2} + \tfrac{\varepsilon}{2} =\varepsilon.
\end{align*}
Then $\mathcal{S}: \mathcal{O}_{\varepsilon} \to \mathcal{O}_{\varepsilon}$ is well defined.

\justify

Let us prove that $\mc{S}$ is a compact map. 
Let $(u_{n}) \subset \mc{O}_{\varepsilon} \subset H^{s}_{G}(M)$ be a bounded sequence; we know that $\mc{S}(u_{n})$ is bounded in $H^{s}(M)$, but since the above computation implies $\|\mathcal{S}(u)\|_{H^{s}(M)} \leq \|\mathcal{S}(u)\|_{H^{r}(M)} \leq  \varepsilon$, it follows that $\mc{S}(u_{n})$ is a bounded sequence in $H^{r}_{G}(M)$. By the chain of embeddings \eqref{teorema final tesis embedding 4}, the embedding $H^{r}_{G}(M) \hookrightarrow H^{s}_{G}(M)$ is compact; therefore the sequence $\iota(\mc{S}(u_{n}))$ admits a subsequence $\iota(\mc{S}(u_{n_{k}}))$ strongly convergent in $H^{s}_{G}(M)$.
Hence $\mc{S}(u_{n_{k}})$ is a subsequence that converges strongly in $H^{s}_{G}(M)$; therefore $\mc{S}$ is a compact map.
If we further prove that $\mc{S}$ is continuous, we can conclude by Schauder's Theorem \cite[c.11.2]{Trudinger} that the operator $\mc{S}$ has a fixed point, which implies that there exists a solution to the system $L_{c}u=\mc{A}(\cdot,u)$.  
Let us prove that $\mc{S}$ is a continuous map.
Let $(u_{n}) \subset \mathcal{O}_{\varepsilon}$ with $u_{n} \to u$ in $H^{s}_{G}(M)$. Using the linearity of the system $L_{c}u=\mathcal{A}(\cdot,u)$. 
By the embeddings \eqref{teorema final tesis embedding 4} and \eqref{teorema final tesis embedding 5}, we have the following computations
\begin{align*}
   \norm{\mathcal{S}(u_{n}) - \mathcal{S}(u) }^{2}_{H^{s}(M)} &= \norm{\widetilde{u}_{n} - \widetilde{u} }_{H^{s}(M)}^{2} \leq \norm{\widetilde{u}_{n} - \widetilde{u}}_{\widetilde{H}_{f}(M)}^{2} = \norm{\mathcal{A}(\cdot,u_{n}) - \mathcal{A}(\cdot,u)}^{2}_{L^{2}(M)} \\
    &= \int_{M}|\mathcal{A}(x,u_{n})  - \mathcal{A}(x,u)|^{2}dV_{g}=\int_{M}\left| \int_{0}^{1}\frac{d}{dt}\mathcal{A}(x,tu_{n} + (1-t)u )dt  \right|^{2}dV_{g}\\\
    &= \int_{M}|u_{n}(x)-u(x)|^{2}\left|  \int_{0}^{1}\frac{\partial}{\partial y}\mathcal{A}(x,tu_{n} + (1-t)u )dt    \right|^{2}dV_{g}
\end{align*}
We define the function $k_{x}: [0,1] \to \mathbb{R}, k_{x}(t)=\frac{\partial}{\partial y}\mathcal{A}(x, tu_{n}(x) 
   +  (1-t)u(x))$ with $x \in M$, since $U(\cdot, u)$ is $C^{1}$ in the second variable, $k_{x}(t)$ is a well defined Carathéodory function.
   Note that
\begin{equation*}
\Big|\int_{0}^{1}k_{x}(t)dt \Big| \leq \int_{0}^{1} |k_{x}(t)|dt \leq \int_{0}^{1}\left( \sup_{t \in [0,1]}|k_{x}(t)|  \right)dt = \sup_{t \in [0,1]}|k_{x}(t)|
\end{equation*}
Since $k_{x}(t)$ is continuous in a compact set, the maximum is attained for each $x \in M$, i.e., for each fixed $x_{0} \in M$,
there exists $t_{0} \in [0,1]$ such that $k_{x_{0}}$ attains the supremum in $t_{0}$. In fact there exists a Borel measurable function $t(\cdot): M \to [0,1]$ by \cite[t.18.13]{AliprantisBorder}
such that for each $x^{\prime} \in M$, $t(x^{\prime})$ 
satisfies $\sup_{t \in [0,1]}|k_{x^{\prime}}(t)|=|k_{x^{\prime}}(t(x^{\prime}))|$. By the inequalities (\ref{desigualdad de crecimiento}), $0 \leq t(x) \leq 1$, Hölder inequality and the embedding \eqref{teorema final tesis embedding 6} we get for some constants $C, \widetilde{C}> 0$ that
\begin{align*}
  \norm{\mathcal{S}(u_{n}) - \mathcal{S}(u)}^{2}_{H^{s}(M)} 
  &\leq  \int_{M}|u_{n}(x) - u(x)|^{2}\left\{\sup_{t \in [0,1]}\left| \frac{\partial}{\partial y}\mathcal{A}(x,tu_{n}(x) + (1-t)u(x))\right|     \right\}^{2}dV_{g} \\ 
  &\leq  \int_{M}|u_{n}(x) - u(x)|^{2}\left\{\left| \frac{\partial}{\partial y}\mathcal{A}(x,t(x)u_{n}(x) + (1-t(x))u(x))\right|     \right\}^{2}dV_{g} \\ 
  &\leq C\int_{M}|u_{n}(x) - u(x)|^{2} \left( 1 + |t(x)u_{n}(x) + (1-t(x))u(x)|^{\beta} \right)^{2}dV_{g}\\
  &\leq \widetilde{C}\int_{M}|u_{n}(x) - u(x)|^{2}\left(  1 + |u_{n}(x)|^{2 \beta} + |u(x)|^{2 \beta}  \right)dV_{g}\\
&\leq \widetilde{C} \left( \norm{u_n - u}_{L^2(M)}^2 + \left(\|u_{n}\|_{L^{2(\beta+1)}(M)}^{2\beta} + \|u\|_{L^{2(\beta+1)}(M)}^{2\beta}\right)\norm{u_{n} - u}_{L^{2(\beta+1)}(M)}^{2}\right)\\
 &\leq \widetilde{C}\left( \left(1+2\varepsilon^{2\beta} \right)\norm{u_{n} - u}_{H^{s}(M)}^{2} \right). 
\end{align*}
Hence $\widetilde{u}_n \to \widetilde{u}$ in $H^{s}(M)$, so $\mathcal{S}$ is continuous. By Schauder's fixed point Theorem \cite[c.11.2]{Trudinger} there exists a fixed point $\phi \in \mathcal{O}_\varepsilon \subset H^s_G(M)$ such that $\phi = \mathcal{S}(\phi) \in \widetilde{H}_{f,G}(M,g)$, which is a non-trivial solution of \eqref{sistema lineal caso Riemaniano con accion}.   
\end{proof}

\newpage

\section*{Acknowledgements} I would like to express my gratitude to my advisor, Dr. Enrique G. Reyes, and to Universidad de Santiago de Chile for their financial support and ANID through Fondecyt Grant No. 1241719.

\end{document}